\documentclass[letterpaper,journal]{article}
\usepackage{algorithmic}
\usepackage{algorithm}
\usepackage{array}
\usepackage[caption=false,font=normalsize,labelfont=sf,textfont=sf]{subfig}
\usepackage{textcomp}
\usepackage{stfloats}
\usepackage{url}
\usepackage{verbatim}
\usepackage{graphicx}
\usepackage{cite}
\usepackage{amsmath,amssymb,amsfonts}
\usepackage{amsthm}
\usepackage{algorithmic}
\usepackage{graphicx}
\usepackage{algorithm,algorithmic}
\usepackage{hyperref}
\hypersetup{hidelinks}
\usepackage{textcomp}
\newtheorem{definition}{Definition}
\usepackage{mathtools}
\DeclarePairedDelimiter\ceil{\lceil}{\rceil}
\newtheorem{theorem}{Theorem}
\newtheorem{remark}{Remark}
\newtheorem{proposition}{Proposition}
\newtheorem{assumption}{Assumption}
\usepackage{siunitx}

\renewcommand{\b}[1]{{\boldsymbol{#1}}}
\def\BibTeX{{\rm B\kern-.05em{\sc i\kern-.025em b}\kern-.08em
    T\kern-.1667em\lower.7ex\hbox{E}\kern-.125emX}}
\begin{document}

\title{ROBBO: An Efficient Method for Pareto Front Estimation with Guaranteed Accuracy}
\author{Roberto Boffadossi, Marco Leonesio, Lorenzo Fagiano
\thanks{This research has been supported by the EU project E2COMATION (H2020 - Agreement N. 958410), by the European Union–Next Generation EU in the context of the project PNRR M4C2, Investimento 1.3 DD. 341 del 15 marzo 2022 – NEST – Network 4 Energy Sustainable Transition – Spoke 2 - PE00000021 - D43C22003090001, by Fondazione Cariplo under grant n. 2022-2005, project ``NextWind - Advanced control solutions for large scale Airborne Wind Energy Systems'', and by the Italian Ministry of University and Research through the European Union – NextGenerationEU fund, project P2022927H7
``DeepAirborne – Advanced Modeling, Control and Design Optimization Methods for Deep Offshore Airborne Wind Energy''.\\ \textit{(Corresponding author: R. Boffadossi, L. Fagiano)}}
\thanks{R. Boffadossi and M. Leonesio are with 
the Institute of Intelligent Industrial Technologies and Systems for Advanced Manufacturing, National Research Council, Milan, 20133 Italy (e-mail: roberto.boffadossi@stiima.cnr.it; marco.leonesio@stiima.cnr.it) and with the Department of Electronics, Information and Bioengineering, Politecnico di Milano, Milan, 20133 Italy.}
\thanks{L. Fagiano is with the Department of Electronics, Information and Bioengineering, Politecnico di Milano, Milan, 20133 Italy (e-mail: lorenzo.fagiano@polimi.it)}}

\markboth{\hskip20pc}
{Boffadossi \MakeLowercase{\textit{et al.}}: ROBBO: An Efficient Method for Pareto Front Estimation with Guaranteed Accuracy}

\maketitle

\begin{abstract}
A new method to estimate the Pareto Front (PF) in bi-objective optimization problems is presented. Assuming a continuous PF, the approach, named ROBBO (RObust and Balanced Bi-objective Optimization), needs to sample at most a finite, pre-computed number of PF points. Upon termination, it guarantees that the worst-case approximation error lies within a desired tolerance range, predefined by the decision maker, for each of the two objective functions. Theoretical results are derived, about the worst-case number of PF samples required to guarantee the wanted accuracy, both in general and for specific sampling methods from the literature. A comparative analysis, both theoretical and numerical, demonstrates the superiority of the proposed method with respect to popular ones. The approach is finally showcased in a constrained path-following problem for a 2-axis positioning system and in a steady-state optimization problem for a Continuous-flow Stirred Tank Reactor. An open demo implementation of ROBBO is made available online.
\end{abstract}

\section{Introduction}
Bi-objective optimization problems (BOPs) are ubiquitous in science, engineering, and economics \cite{BOP1,BOP3,BOP4}. 
In control, BOPs are extremely common, because control tuning yields a wide range of possible trade-offs. BOPs account for over 50\% of the tuning problems, while about 30\% pertain to three objectives, and 20\% more than three, see \cite{MultiobjectiveMetaheuristicOptimization2020a}.
A typical approach to deal with a BOP is to approximate the set of all the non-dominated feasible points in the objectives (or criteria) space, named Pareto Front (PF), and let the Decision Maker (DM) manually select the desired one, based on various application-specific considerations. Such an \textit{a posteriori} articulation of preferences \cite{marlerSurveyMultiobjectiveOptimization2004} is the most common and relevant in real world settings, and also the most documented case for control problems, as shown in \cite{MultiobjectiveMetaheuristicOptimization2020a}.
\noindent There are three categories of  approximation strategies: pointwise, piecewise, and others.

The first class corresponds to computing a finite ``solution set'' that acts as a surrogate of the whole PF. 
The most common methods provide a point-cloud approximation by sampling the PF, each sample being the solution of an auxiliary single-objective optimization problem with scalarization techniques to find a specific trade-off (see, for instance,  \cite{dasCloserLookDrawbacks1997,dasNormalBoundaryIntersectionNew1998,ismail-yahayaEffectiveGenerationPareto2002, CHIANDUSSI2012912,alianofilhoExactScalarizationMethod2021}).
Evolutionary strategies or meta-heuristics have been proposed as well (\cite{emmerichTutorialMultiobjectiveOptimization2018}; \cite{pereiraReviewMultiobjectiveOptimization2022a}), proving to be effective even for irregular PFs. However, the complete set of Pareto optimal solutions can involve an infinite number of potential trade-offs, thus making it impossible to be fully described by a finite set of points. Moreover, reaching a good approximation is impractical when the sampling process is expensive, for example requiring time-consuming simulations, global optimization, or experiments.
The second class, of piecewise approximation methods, is a direct extension of the previous one, consisting of different interpolation/fitting strategies of the elements of the solution set: piecewise linear (\cite{cohonGeneratingMultiobjectiveTradeoffs1979}; \cite{kimAdaptiveWeightedSum2006}), outer,  inner, and ``sandwich'' approximations (\cite{klamrothUnbiasedApproximationMulticriteria2003}; \cite{bokrantzAlgorithmApproximatingConvex2013a}), quadratic or cubic \cite{ruzikaApproximationMethodsMultiobjective2005}. To refine the approximation, these methods exploit iterative techniques to sample the PF considering a chosen error metric. 

The last category includes meta-models or higher-order functions to provide a surrogate PF. Examples are enclosure or box coverage, like  \cite{eichfelderAdvancementsComputationEnclosures2023b,eichfelderApproximationAlgorithmMultiobjective2022,
payneInteractiveRectangleElimination1980, payne1993efficient}, which approximate the PF as the union of rectangles in between the sampled PF points and evaluate the approximation quality by the maximum size of such rectangles, or Gaussian Processes like \cite{campigottoActiveLearningPareto2014,binoisQuantifyingUncertaintyPareto2015}, which exploit a sampling strategy based on the estimation of uncertainty derived by strategically collecting new PF observations, or combination of local manifolds \cite{tianLocalModelBasedPareto2023a}, which has shown to be particularly effective when dealing with complex, irregular PFs.

Finding an accurate estimate of the PF with as little computational or experimental effort as possible, i.e. the smallest number of samples, is important for an effective and efficient decision-making process. 
At the same time, providing a guarantee on the estimation error is crucial to make an informed decision. 
Ideally, the DM would like to have a method that a) builds an estimated PF with proven guaranteed error with respect to the real one, with b) the possibility to specify the error tolerances for each criterion, c) knowing early on the worst-case number of PF samples needed to guarantee the desired accuracy, and d) paired with a ``realization mechanism'' that, once a point on the estimated PF has been selected, returns a point on the real PF that satisfies the tolerance requirements. Notwithstanding the large number of PF estimation approaches in the literature, to the best of our knowledge, none of them addresses all these requirements. We refer to such a problem as a \textit{robust and balanced} PF estimation-realization: robust, because it provides accuracy guarantees, and balanced, because it constrains the ratio of the estimation errors of the two objectives to a desired value.

The main contribution of this paper is to provide such a method. First, we introduce a linear, invertible transformation of the PF that casts the problem as that of approximating a conveniently chosen univariate scalar function. Such a transformation takes into account the error tolerances defined by the DM. Under a continuity assumption, we derive the set of all possible PFs (named Feasible Pareto Function Set) that are compatible with the sampled points, as well as its tight upper and lower bounds. A first theoretical result establishes the necessary and sufficient conditions on the PF samples to guarantee the wanted accuracy, for any possible approximation belonging to the feasible set. Then, our new method is introduced, named ROBBO (RObust and Balanced Bi-objective Optimization). It is an iterative algorithm that uses the point-wise distance between the bounds to select the most informative next sample, among a finite set of suitably chosen candidates, thus obtaining an efficient algorithm where the worst-case uncertainty is quickly reduced to below the desired threshold. 
After presenting the method and deriving its theoretical properties, we compare it theoretically and numerically with popular alternative methods, and test it on a constrained path-following problem for a 2-axis positioning system, and on the computation of the economically optimal steady state of a Continuous-flow Stirred Tank Reactor (CSTR).

The paper is organized as follows. Section \ref{S:problem_formulation} includes preliminary notions and problem formulation. Sections \ref{S:theoretical conditions}-
\ref{sec:extensions, comparisons} present the theoretical results,  proposed approach,  extensions and comparisons. Section \ref{S:Results} contains the simulation results and Section \ref{S:conclusions} concludes the paper.

\textbf{\textit{Notation}}. A bold-faced symbol like $\b{z}$ denotes a vector, $z_i$ its $i^\text{th}$ entry, and $\b{z}^T$ its transpose.  
$\mathbb{R}^2_{\succeq 0}$ and $\mathbb{R}^2_{\preceq 0}$ are the set of vectors in $\mathbb{R}^2$ with non-negative and with non-positive entries. The \textit{Pareto Cone} is $\mathbb{R}^2_{\succ0}\doteq\mathbb{R}^2_{\succeq 0}\setminus\{[0,0]^T\}$. Similarly, we define $\mathbb{R}^2_{\prec0}\doteq\mathbb{R}^2_{\preceq 0}\setminus\{[0,0]^T\}$. The Minkowski sum of two sets $A$ and $B$ is $A\oplus B\doteq\{a+b\;|\;a\in A \wedge b\in B\}$.
\section{Problem statement}\label{S:problem_formulation}
We consider the BOP:
\begin{equation}\label{MOP}
    \min\limits_{\b{x}\in{S}}\boldsymbol{f}(\b{x})
\end{equation}
where function $\b{f}(\b{x}):\mathbb{R}^n\rightarrow\mathbb{R}^2$ returns the two objectives (or criteria) $z_{1,2}=f_{1,2}(\b{x}):\mathbb{R}^n\rightarrow\mathbb{R}$, and
the design variables are restricted to the feasible design set $S\subset{\mathbb{R}^n}$. We denote the feasible criterion set with $Z\doteq \{\b{f}(\b{x}),\, \b{x}\in{S}\}\subset\mathbb{R}^2$. We further denote a specific sample in $S$, with an associated index $j$, as $\b{x}^{(j)}$, and the corresponding objective vector as $\b{z}^{(j)}= \b{f}(\b{x}^{(j)})=[z^{(j)}_1,z^{(j)}_2]^T$.
We recall a few fundamental definitions \cite{emmerichConeBasedHypervolumeIndicators2013}, \cite{nakayamaBasicConceptsMultiobjective2009}.
\begin{definition}\label{dominance}
  A point $\b{z}^{(1)}\in{Z} $ dominates a point $\b{z}^{(2)}\in{Z} $, with notation $\b{z}^{(1)}\prec\b{z}^{(2)}$, iff $\b{z}^{(2)} \in \{\b{z}^{(1)}\}\oplus\mathbb{R}^2_{\succ0} $. That is, $\b{z}^{(1)}$ presents better performance in one objective and at least equal in the other.
\end{definition}
\begin{definition}
    The set of Pareto optimal decisions is $X^*\doteq\{\b{x}^*\in{S}\,|\,\nexists \b{x}\in S, \b{x}\neq\b{x}^*:\;\b{f}(\b{x})\prec\b{f}(\b{x}^*)\}$. 
\end{definition}        
\begin{definition}
        The Pareto Front (PF) is the set of all non-dominated points in $Z$,
        $\mathcal{P}\doteq\{\b{z}^*=\b{f}(\b{x}^*)\;:\;\b{x}^*\in X^*\}$.
\end{definition}
We refer to $\b{x}^*$ and $\b{z}^*$ as a Pareto optimal point in the decision space and in the criterion space, respectively.\\
\noindent We next introduce two specific points belonging to the PF, named \textit{anchor points},
$\b{z}^{*(a1)}$ and
$\b{z}^{*(a2)}$,
obtained by solving the following auxiliary problems:
\begin{subequations}\label{eq:anchor points}
\begin{gather}
        \b{z}^{*(a1)}= \b{f}(\b{x}^{*(a1)}): \b{x}^{*(a1)}=\arg\min_{\b{x}\in S}f_1(\b{x})\label{eq:anchor points-1}\\
        \b{z}^{*(a2)}= \b{f}(\b{x}^{*(a2)}): \b{x}^{*(a2)}=\arg\min_{\b{x}\in S}f_2(\b{x})\label{eq:anchor points-2}
\end{gather}
\end{subequations}
\begin{assumption}\label{ass.PF}
    $\mathcal{P}$ is a continuous and compact curve in $\mathbb{R}^2$.
\end{assumption}

\begin{figure}
\centerline{\includegraphics[trim=353 280 280 27,clip,scale=0.78]{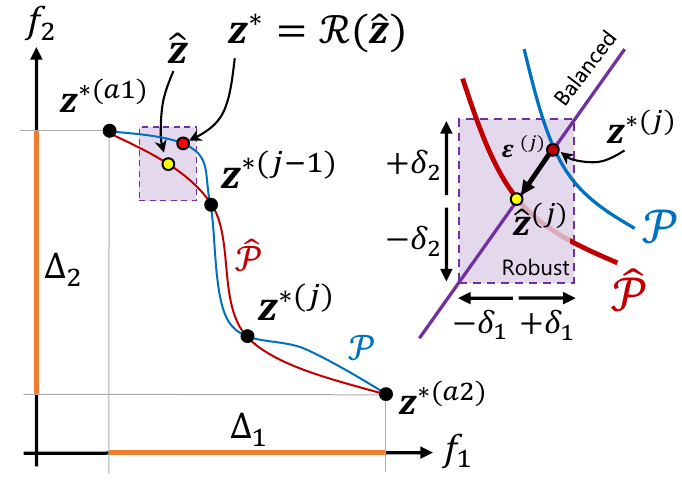}}
\caption{Robust and balanced PF estimation-realization: Pareto front (blue line) and its approximation (red line), candidate selection (yellow point) and its realization (red point) within the tolerance region (dashed violet box).}
\label{fig.setup}
\end{figure}

\noindent Satisfaction of Assumption \ref{ass.PF} depends on the application at hand and on the modeling choices made when defining $\boldsymbol{f}$ and $S$. Compactness is a very mild condition, since it holds whenever the anchor points \eqref{eq:anchor points} are finite, which is always the case in sensible engineering applications. Note that compactness of the PF does not require compactness of the feasible design set $S$, nor of the feasible criterion set $Z$.
Also, continuity is a mild assumption in real-world engineering problems, where discontinuous PFs are not frequently encountered  \cite{logistFastParetoSet2010a}. Further, we discuss possible strategies to handle the presence of discontinuities in our approach
at the end of Section \ref{S:Sampling}. \\ 
The entries of $\b{z}^{*(a1)},\,\b{z}^{*(a2)}$ form the upper and lower extremes of the closed and bounded intervals obtained by projecting $\mathcal{P}$ on the Cartesian axes in $\mathbb{R}^2$. We indicate the lengths of such intervals with $\Delta_1\doteq z^{*(a2)}_1-z^{*(a1)}_1$ and $\Delta_2\doteq z^{*(a1)}_2-z^{*(a2)}_2$, see Fig. \ref{fig.setup}.

We now consider an estimate $\hat{\mathcal{P}}\approx\mathcal{P}$ (approximated PF). We denote a point in $\hat{\mathcal{P}}$ with $\hat{\b{z}}$, and introduce a concept named \textit{realization}. This is a mapping $\mathcal{R}:\hat{\mathcal{P}}\rightarrow\mathcal{P}$, i.e., for any $\hat{\b{z}}^{(j)}\in\hat{\mathcal{P}}$ we have that $\mathcal{R}(\hat{\b{z}}^{(j)})\in\mathcal{P}$. In a decision making perspective we will refer to $\hat{\b{z}}^{(j)}$ as the \textit{candidate}, i.e. the point selected on the approximation $\hat{\mathcal{P}}$ by the DM, and to $\b{z}^{*(j)}=\mathcal{R}(\hat{\b{z}}^{(j)})$ as the \textit{realization of the candidate}, i.e. the related point on the actual PF.
We denote the realization error as the vector
$\b{\varepsilon}=\b{\hat{z}}-\b{z}^{*}.$

We consider a scenario where the DM wants to derive $\hat{\mathcal{P}}$ and $\mathcal{R}$ such that the following property of the realization error holds:
\begin{equation}\label{eq:robust requirement}
\forall\hat{\b{z}}\in\hat{\mathcal{P}},\,\b{z}^{*}=\mathcal{R}(\hat{\b{z}}),\,\begin{array}{l}
|\varepsilon_1|<\delta_1\wedge |\varepsilon_2|<\delta_2
\end{array}
\end{equation}
where $\delta_1,\,\delta_2$ are defined by the DM according to the wanted accuracy/tolerance for each criterion. As anticipated in the Introduction, we refer to a pair $(\hat{\mathcal{P}},\,\mathcal{R})$ satisfying \eqref{eq:robust requirement} as a \textit{robust} PF estimate-realization.

We further consider that the approximated PF is obtained on the basis of a finite number $M$ of PF samples, which include to a minimum the anchor points (thus, $M\geq 2$), forming the data-set  $\mathcal{D}=\left\{\b{z}^{*(a1)},\b{z}^{*(1)},\ldots,\b{z}^{*(M-2)},\b{z}^{*(a2)}\right\}$. Without loss of generality, we consider the elements $\b{z}^{*(j)}\in\mathcal{D}$ to be ranked according to the $f_1$-coordinate so that $z_1^{*(j)}<z_1^{*(j+1)}$. A graphical illustration of the considered setup is presented in Fig. \ref{fig.setup}.

Since obtaining each sample in $\mathcal{D}$ may require the solution of complex and time consuming simulations, global optimization and/or experiments, it is of interest to the DM to have an estimate-realization that is not only robust, but also requires a ``low'' value of $M$.
From the theoretical standpoint, our goal is thus to investigate what are the conditions on $\mathcal{D}$ to satisfy \eqref{eq:robust requirement}, whether there is a minimum number $M$ of suitably chosen samples such that  \eqref{eq:robust requirement} is surely satisfied, and whether there are estimation-realization algorithms that are particularly \textit{efficient} in terms of worst-case value of $M$, i.e. they require a low number of samples to provide the accuracy guarantees.

As we show in Section \ref{S:theoretical conditions}, we can cast the PF estimation problem into a function approximation one. Then, given a data-set $\mathcal{D}$, one can derive the set of functions that are compatible with the observed PF samples. We name this set the Feasible Pareto Function Set (${FPFS}_{\mathcal{D}}$). With this concept in mind, we formulate the problems \textbf{P1}-\textbf{P3} addressed in this paper as follows:

\textbf{P1}. Derive necessary and sufficient conditions on $\mathcal{D}$ such that, for any PF $\mathcal{P}$ and any estimate $\hat{\mathcal{P}}$ belonging to the corresponding ${FPFS}_{\mathcal{D}}$, there exists a realization $\mathcal{R}$ such that \eqref{eq:robust requirement} is satisfied.

\textbf{P2}. Derive the smallest number $\underline{M}$ of PF samples required to guarantee the conditions of \textbf{P1}.

\textbf{P3}. Obtain an algorithm that builds {a specific estimate-realization pair} $(\hat{\mathcal{P}},\mathcal{R})$ satisfying \eqref{eq:robust requirement} with a finite number of samples $M<\underline{M}$, which is the smallest possible in the worst case.

We remark that problems \textbf{P1}-\textbf{P2} consider any possible (admissible, i.e., belonging to the ${FPFS}_{\mathcal{D}}$) estimate $\hat{\mathcal{P}}$ and, since the PF is not known except for its samples, any possible (admissible) PF. The value $\underline{M}$ shall thus be valid for all possible combinations of $\hat{\mathcal{P}},\mathcal{P}\in {FPFS}_{\mathcal{D}}$.\\
On the other hand, in problem \textbf{P3} we aim for a specific estimate $\hat{\mathcal{P}}$, thus the ``worst case'' scenario mentioned in \textbf{P3} refers to the uncertainty about $\mathcal{P}$ only. This is why the value of $M$ in \textbf{P3} is generally smaller than $\underline{M}$, and we look for the estimation approach that minimizes it.
We provide the solutions to problems \textbf{P1}-\textbf{P3} in the next sections.

\section{Robust PF estimation - theoretical results} \label{S:theoretical conditions}
\subsection{Feasible PF Set and its optimal bounds}
\begin{figure}[!bht]
\centerline{\includegraphics[trim=100 250 420 50,clip,scale=0.565]{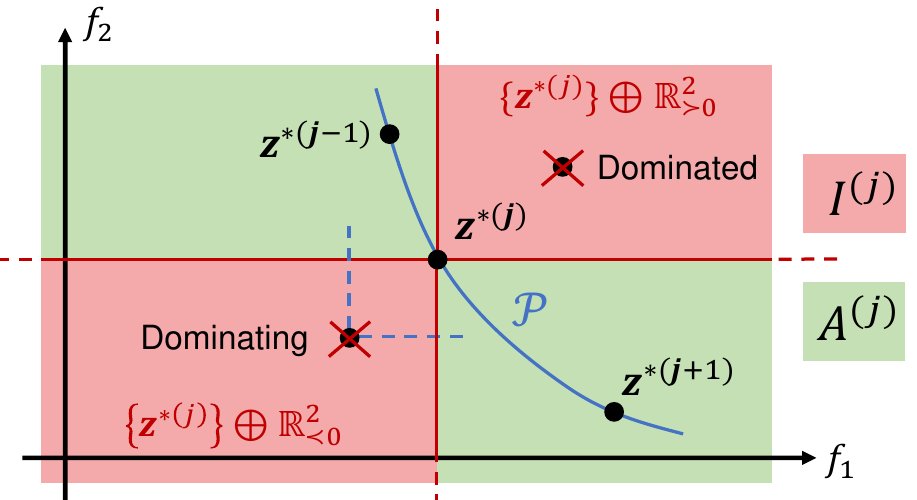}}
\caption{Subsets ${I}^{(j)}$ (red) and ${A}^{(j)}$ (green) pertaining to a point $\b{z}^{*(j)}$}
\label{fig.zones}
\end{figure}

According to the Pareto dominance relation, each PF point $\b{z}^{*(j)}\in\mathcal{D}$ defines the following two subsets:
 \begin{equation}\label{eq:inactive-active sets}
     I^{(j)} \doteq \{\b{z}^{*(j)}\}\oplus\mathbb{R}^2_{\succ0}\;\cup\;\{\b{z}^{*(j)}\}\oplus\mathbb{R}^2_{\prec0},\;\;A^{(j)} \doteq\mathbb{R}^2 \setminus I^{(j)}
 \end{equation}
The set $I^{(j)}$ is the subset of the criterion space in which no points of the PF can be present, because they would either be dominated or dominate point $\b{z}^{*(j)}$, contradicting the fact that the latter belongs to the PF, see Fig. \ref{fig.zones}.
Thus, $I^{(j)}\cap \mathcal{P} = \emptyset,\,\forall j$.
Regarding the anchor points, the corresponding sets  $I^{(a1)},\,I^{(a2)}$ are larger, i.e.:
 \begin{equation}\nonumber
 \begin{array}{rcl}
     I^{(a1)} &\doteq& \{\b{z}^{*(a1)}\}\oplus\mathbb{R}^2_{\succ0}\;\cup\;\{\b{z}^{*(a1)}\}\oplus\mathbb{R}^2_{\prec0}\\
     &&\cup\;\{\b{z}:z_1< z_1^{*(a1)}\}  \\
     I^{(a2)} &\doteq& \{\b{z}^{*(a2)}\}\oplus\mathbb{R}^2_{\succ0}\;\cup\;\{\b{z}^{*(a2)}\}\oplus\mathbb{R}^2_{\prec0}\\
     &&\cup\;\{\b{z}:z_2< z_2^{*(a2)}\}
\end{array}
 \end{equation}
 while $A^{(a1)},A^{(a2)}$ are still their complements, as in \eqref{eq:inactive-active sets}.

Considering now the set $\mathcal{D}$ altogether, we can partition the objective domain in the set $\mathcal{I}_{\mathcal{D}}\subset\mathbb{R}^2$ and its complement, $\mathcal{A}_{\mathcal{D}}$:
\begin{equation}\label{eq:inactive-active zones global}
        \mathcal{I}_{\mathcal{D}}\doteq\bigcup\limits_{\b{z}^{*(j)}\in\mathcal{D}}I^{(j)},\;\;\;
        \mathcal{A}_{\mathcal{D}}\doteq\mathbb{R}^2\setminus\mathcal{I}_{\mathcal{D}}
\end{equation}
As presented in Fig. \ref{fig.rotation} (top), $\mathcal{A}_{\mathcal{D}}$ is the subset of the criterion space that is guaranteed to contained the actual PF, according to the information provided by the data $\mathcal{D}$. It also corresponds to the intersection
$\mathcal{A}_{\mathcal{D}}=\bigcap\limits_{\b{z}^{*(j)}\in\mathcal{D}}A^{(j)}$. 
By construction, we have that $\mathcal{P}\in\mathcal{A}_{\mathcal{D}}$ and $\mathcal{I}_{\mathcal{D}}\cap\mathcal{P}=\emptyset$.\\
The set $\mathcal{A}_{\mathcal{D}}$ is instrumental to define what we call the Feasible Pareto Function Set,  ${FPFS}_{\mathcal{D}}$.
We first convert the original coordinate system $(f_1,f_2)$ into a new one, denoted by $(v,q)$, by the following transformation matrix $W$:
\begin{equation}\label{eq:transformation}
    \begin{bmatrix}
v\\
q
\end{bmatrix}
=\underbrace{\frac{\sqrt{2}}{2}
\begin{bmatrix}
1 & -1 \\
1 & 1\\ \end{bmatrix}
\begin{bmatrix}
\frac{1}{\delta_1} & 0 \\
0 & \frac{1}{\delta_2}\\ \end{bmatrix}}\limits_{W}
\begin{bmatrix}
f_1\\
f_2
\end{bmatrix}
\end{equation}
Namely, the transformation consists of normalizing the two objective functions by the corresponding desired error tolerances, and of rotating the scaled axes so that $f_1/\delta_1$ overlaps with the bisector of the first quadrant of the new coordinates, see Fig. \ref{fig.rotation}.
This transformation allows us to study the problem irrespective of the specific error bounds $\delta_1,\delta_2$. 
For any PF point $[z^*_1,z^*_2]^T$, we denote the corresponding point in the transformed coordinates with
\[\left[\begin{array}{c}
     z^*_v  \\
     z^*_q 
\end{array}\right]=W
\left[\begin{array}{c}
     z^*_1  \\
     z^*_2 
\end{array}\right]
\]
The transformed anchor points are thus $[z_v^{*(a1)},z_q^{*(a1)}]^T=W\b{z}^{*(a1)}$ and $[z_v^{*(a2)},z_q^{*(a2)}]^T=W\b{z}^{*(a2)}$. Since the transformation corresponds to a  scaling with positive coefficients and a clockwise rotation by $\frac{\pi}{4}$ of the objectives' axes $(f_1,f_2)$, and recalling that the points in $\mathcal{D}$ are ordered by increasing $f_1$-values, note that the points in the rotated data-set are still ordered by increasing values of their first coordinate, i.e. $z_v^{*(a1)}<z_v^{*(1)}<\ldots<z_v^{*(M-2)}<z_v^{*(a2)}$.\\
Moreover, note that, by the Pareto dominance relation, in general for any two points $\b{z}^{*(j)},\,\b{z}^{*(\ell)}\in\mathcal{P}$  we have $z_1^{*(j)}<z_1^{*(\ell)} \Leftrightarrow z_2^{*(j)}>z_2^{*(\ell)}$. Thus, $\mathcal{P}$ can be seen as a strictly monotonically decreasing function relating $z_2^{*}$ with $z_1^{*}$. Since in \eqref{eq:transformation} we apply a linear and invertible transformation, we can still interpret the transformed front as a function, denoted by $h$, such that $z^*_q=h(z^*_v)$. Under Assumption \ref{ass.PF}, $h$ is continuous (i.e., $h\in C^0$). Its domain is given by the interval delimited by the transformed anchor points projected on the $v$-axis, i.e. $[z_v^{*(a1)},z_v^{*(a2)}]$. Thus, the actual Pareto Front $\mathcal{P}$ can be also written as:
\begin{equation}\label{eq:PF as a function}
\mathcal{P}=\left\{W^{-1}[v,h(v)]^T,\,\forall v\in [z_v^{*(a1)},z_v^{*(a2)}]\right\}
\end{equation}
which for simplicity we denote in short as $\mathcal{P}=W^{-1}h$. 
We are now in position to define the ${FPFS}_{\mathcal{D}}$, i.e., the set of all continuous functions $h$ (i.e., all PF $\mathcal{P}$ according to \eqref{eq:PF as a function}) compatible with the samples in $\mathcal{D}$:
\begin{figure}
\centerline{\includegraphics[trim=310 30 315 15,clip,scale=0.77]{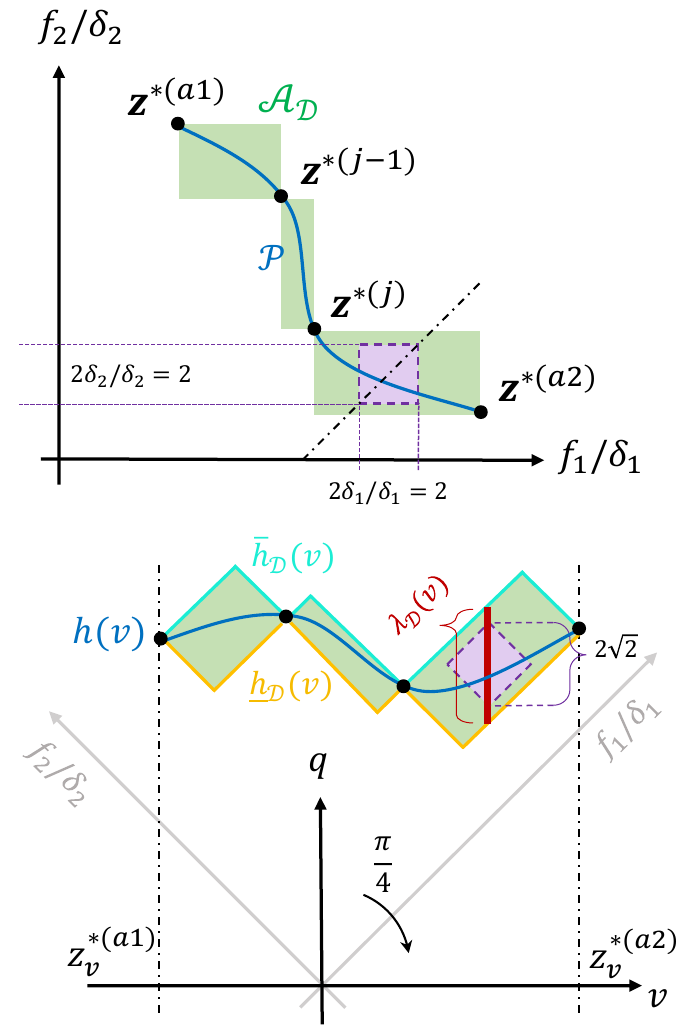}}
\caption{Top figure: Pareto front as a curve $\mathcal{P}$ in the normalized coordinates $(f_1/\delta_1,f_2/\delta_2)$. Bottom figure: Pareto front as a function $h(v)$ in $(v,q)$ coordinates. The dashed violet boxes represent the tolerance region before and after the rotation}
\label{fig.rotation}
\end{figure}
\begin{equation*}
\begin{array}{l}
{FPFS}_{\mathcal{D}}\doteq\\\left\{h\in C^0,h:[z_v^{*(a1)},z_v^{*(a2)}]\rightarrow \mathbb{R}\;\vert\;W^{-1}h\in\mathcal{A}_{\mathcal{D}}\right\}
\end{array}
\end{equation*}
\noindent In this framework, estimating the PF is thus equivalent to deriving a function $\hat{h}\approx h$ in the rotated coordinates. For a given function $\hat{h}$, one can in fact obtain the corresponding estimated PF as $\hat{\mathcal{P}}=W^{-1}\hat{h}$.\\
With a slight abuse of notation, we indicate with $\mathcal{P},\hat{\mathcal{P}}\in {FPFS}_{\mathcal{D}}$ a PF (or estimated PF) whose corresponding function $h,\hat{h} \in {FPFS}_{\mathcal{D}}$.
Considering the function set ${FPFS}_{\mathcal{D}}$ allows us to compute tight bounds on the estimation error. We do this in a Set Membership function approximation framework \cite{MILANESE2004957}, by first deriving the so-called optimal upper and lower bounds to the ${FPFS}_{\mathcal{D}}$, $\bar{h}_{\mathcal{D}}$  and $\underline{h}_{\mathcal{D}}$, defined as: 
\begin{equation*}
    \begin{split}
        &\bar{h}_{\mathcal{D}}:[z_v^{*(a1)},z_v^{*(a2)}]\rightarrow \mathbb{R},\,\bar{h}_{\mathcal{D}}(v)\doteq\sup_{h\in{{FPFS}_{\mathcal{D}}}}{h(v)}\\
        &\underline{h}_{\mathcal{D}}:[z_v^{*(a1)},z_v^{*(a2)}]\rightarrow \mathbb{R},\,\underline{h}_{\mathcal{D}}(v)\doteq\inf_{h\in{{FPFS}_{\mathcal{D}}}}{h(v)}
    \end{split}
\end{equation*}
These functions are the tightest upper and lower bounds to the values that any function in $FPFS_\mathcal{D}$ can take. For a generic $\bar{v}\in[z_v^{*(a1)},z_v^{*(a2)}]$, denote:
\begin{equation}\label{eq:indexes generating bounds}
\begin{array}{rcl}
\b{z}^{*(j^-)}(\bar{v})&=&\arg\min\limits_{{\b{z}}^{*(j)} \in \mathcal{D}} |\bar{v}-z_v^{*(j)}|\;\text{s.t.}\;z_v^{*(j)}\leq \bar{v}\\
\b{z}^{*(j^+)}(\bar{v})&=&\arg\min\limits_{{\b{z}}^{*(j)} \in \mathcal{D}} |\bar{v}-z_v^{*(j)}|\;\text{s.t.}\;z_v^{*(j)}\geq \bar{v}
\end{array}
\end{equation}
i.e., the samples whose $v$-coordinates represent the extremes of the smallest segment, among those defined by each pair of samples, containing $\bar{v}$. These two samples always exist, since the anchor points themselves are assumed to be part of $\mathcal{D}$, and they coincide only if $\bar{v}=z_v^{*(j)}$ for some $j$. Otherwise, considering the ordering of the data-set by increasing $v$-values, we have that $\b{z}^{*(j^+)}(\bar{v})$ is the point right after $\b{z}^{*(j^-)}(\bar{v})$ in the data set.\\
\begin{theorem}\label{T:bounds}
Given the BOP  \eqref{MOP} and a data-set $\mathcal{D}$, under Assumption \ref{ass.PF} it holds:
\begin{equation}\label{eq:optimal bounds}
\begin{array}{l}
       \bar{h}_{\mathcal{D}}(v)\doteq\min\limits_{\b{z}^{*(j^-)}(v),\,\b{z}^{*(j^+)}(v)}{z_q}^{*(j)}+\vert v-{z_v}^{*(j)}\vert\\
       \underline{h}_{\mathcal{D}}(v)\doteq\max\limits_{\b{z}^{*(j^-)}(v),\,\b{z}^{*(j^+)}(v)}{z_q}^{*(j)}-\vert v-{z_v}^{*(j)}\vert
\end{array}
\end{equation}
\end{theorem}
\begin{proof}
See the Appendix.
\end{proof}
\begin{remark}\label{rem:Lipschitz}
From Theorem \ref{T:bounds}, it is evident that the upper and lower bounds \eqref{eq:optimal bounds} are Lipschitz continuous functions with Lipschitz constant equal to one. Therefore, by a local analysis in the neighborhood of a generic PF point $\b{z}^*$, one can derive that, in the considered problem, $h(v)$ is a Lipschitz continuous function with Lipschitz constant smaller than one (see Fig. \ref{fig.rotation}, bottom graph, for a visualization). This fact could be used to further refine the $FPFS_\mathcal{D}$, by restricting it to such a function class. However, this restriction would not change any of the findings and is thus omitted for simplicity.
\end{remark}

\subsection{Conditions for robustness guarantees}\label{SS:minimum M}
Equipped with the $FPFS_\mathcal{D}$ and its optimal bounds, let us now introduce the local, $\lambda_{\mathcal{D}}(v)$, and global, $\overline{\lambda}_{\mathcal{D}}$, worst-case estimation errors on the basis of the available data $\mathcal{D}$ \cite{MILANESE2004957}:
\begin{subequations}\label{eq:local-global errors}
\begin{gather}
    \lambda_{\mathcal{D}}(v)   \doteq\sup\limits_{\hat{h},\tilde{h}\in {FPFS}_{\mathcal{D}}}|\hat{h}(v)-\tilde{h}(v)|=\bar{h}_{\mathcal{D}}(v)-\underline{h}_{\mathcal{D}}(v)\label{eq:local-global errors-1}\\
    \overline{\lambda}_{\mathcal{D}} \doteq\max\limits_{v\in[z^{*(a1)}_{v},z^{*(a2)}_{v}]}\lambda_{\mathcal{D}}(v) 
    \end{gather}
\end{subequations}
The local worst-case error provides an upper bound to the difference, for a specific $v$, between any two functions in $FPFS_\mathcal{D}$, thus also including the rotated PF $h(v)$, since it belongs to this function set. Equation \eqref{eq:local-global errors-1} can be proven using the definition of the optimal bounds (see e.g. \cite{MILANESE2004957}). The global worst-case error is the maximum local error over all $v$ values in the functions' domain. In the bi-objective optimization problem at hand, since we are dealing with a scalar function of a scalar variable, these quantities can be always computed analytically, using Theorem \ref{T:bounds} and the compactness of the segment $[z^{*(a1)}_{v},z^{*(a2)}_{v}]$. Indeed, computing $\lambda_{\mathcal{D}}(v)$ amounts to applying \eqref{eq:optimal bounds}, while we provide next a result to derive $\overline{\lambda}_{\mathcal{D}}$. 
Let us denote with $R^{(\ell)},\ell=1,\ldots,M-1$, the open intervals on the $v$-axis defined by each pair of consecutive points in $\mathcal{D}$, i.e. $R^{(1)}=(z_v^{*(a1)},z_v^{*(1)}),\,R^{(2)}=(z_v^{*(1)},z_v^{*(2)}),\,\ldots,R^{(M-1)}=(z_v^{*(M-2)},z_v^{*(a2)})$. For a given $\ell=1,\ldots,M-1$, let us denote generically with $\b{z}^{*(\ell-)},\,\b{z}^{*(\ell+)}$ the samples whose $v$-coordinates define the extremes of $R^{(\ell)}$. Finally, let us denote
\[
\begin{array}{rcl}
V^{(\ell)}&=&z_v^{*(\ell+)}-z_v^{*(\ell-)}\\
Q^{(\ell)}&=&z_q^{*(\ell+)}-z_q^{*(\ell-)}
\end{array}
\]
\begin{proposition}\label{prop: max global error} 
Given the BOP problem  \eqref{MOP} and a data-set $\mathcal{D}$, under Assumption \ref{ass.PF} the global worst-case error bound is:
\begin{equation}\label{eq:max global error result}
\overline{\lambda}_{\mathcal{D}}=\max\limits_{\ell=1,\ldots,M-1}V^{(\ell)}-|Q^{(\ell)}|
\end{equation}
\end{proposition}
\begin{proof}
See the Appendix.
\end{proof}

The global bound $\overline{\lambda}_{\mathcal{D}}$
plays a key role in the solution to problem \textbf{P1}, which is provided by the next result. In the remainder, for a generic estimated Pareto point $\hat{\b{z}}$ we will make use of the realization
\begin{equation}\label{eq:SM realization}  
    \b{z}^{*}=\mathcal{R}_{SM}(\hat{\b{z}})\doteq W^{-1}\left[
    \begin{array}{c}
        \hat{z}_v \\
          h(\hat{z}_v)
     \end{array}\right]
     \end{equation}
i.e., such that the $v$-coordinates of the estimated and of the actual Pareto point are the same. The value $h(\hat{z}_v)$ in (\ref{eq:SM realization}) can be computed by solving a suitable auxiliary optimization problem, which we introduce in detail in Section \ref{S:Sampling}.

\begin{theorem}\label{th:2}
    Given the BOP problem  \eqref{MOP}, a data-set $\mathcal{D}$, and a tolerance vector $\boldsymbol{\delta}$, under Assumption \ref{ass.PF} it holds:    
    \begin{equation}\label{eq:main result}
    \begin{array}{c}
        \exists \mathcal{R}:\text{ property \eqref{eq:robust requirement} is satisfied }\forall \hat{\mathcal{P}},\mathcal{P}\in FPFS_\mathcal{D}, \; \\
        \iff\\\overline{\lambda}_{\mathcal{D}}\leq\sqrt{2}
    \end{array}
    \end{equation}    
\end{theorem}
\begin{proof}
See the Appendix.
\end{proof}

We now address problem \textbf{P2}, i.e. to derive the minimum number of samples such that the condition laid out by Theorem \ref{th:2} is met, i.e., $\overline{\lambda}_{\mathcal{D}}\leq\sqrt{2}$.
\noindent
Recall that $\Delta_1=z^{*(a2)}_1-z^{*(a1)}_1$ and $\Delta_2=z^{*(a1)}_2-z^{*(a2)}_2$.
\begin{proposition}\label{prop:minimum M}
    Given the BOP problem  \eqref{MOP} and a tolerance vector $\boldsymbol{\delta}$, under Assumption \ref{ass.PF}, the minimum number of samples required to satisfy condition \eqref{eq:main result} is:
    \begin{equation*}
\underline{M}=\ceil*{\frac{1}{2}\left(\frac{\Delta_1}{\delta_1}+\frac{\Delta_2}{\delta_2}\right)}+1
\end{equation*}
and the corresponding data-set $\mathcal{D}$ shall be evenly allocated in the segment $\left[z^{*(a2)}_v,z^{*(a1)}_v\right]$.
\end{proposition}
\begin{proof}
See the Appendix.
\end{proof}

\begin{remark}\label{rem:uniform points min number} Using Proposition \ref{prop:minimum M}, the DM can evaluate, after computing only the two anchor points, the minimum number of samples that surely (i.e., no matter the shape of the actual PF) attains the robustness conditions of Theorem \ref{th:2}. Note that the anchor points themselves are included in $\underline{M}$, and that the result also prescribes how these points should be sampled, i.e. such that they uniformly partition the interval $[z^{*(a1)}_v,\,z^{*(a2)}_v]$ into segments whose length is at most $\sqrt{2}$.\end{remark}

\section{ROBBO: RObust and Balanced Bi-objective Optimization}\label{S:Approach}
We now address problem \textbf{P3}, i.e., we propose an algorithm, named ROBBO, to obtain an estimate-realization pair $(\hat{\mathcal{P}},\,\mathcal{R})$ that attains the robustness guarantees with a much lower number of samples than what obtained in the general case of Proposition \ref{prop:minimum M}. Moreover, the algorithm attains a balanced error vector, in a sense better specified below.
\subsection{Central estimate, minimum number of samples and balanced error property}
ROBBO adopts the following estimate of the PF for a given set of samples $\mathcal{D}$:
\begin{equation}\label{eq:central_approx}
\begin{array}{rcl}
     \hat{\mathcal{P}}_c&=&W^{-1}\hat{h}_{c}\\
     \hat{h}_{c}&=&\frac{1}{2}({\bar{h}_{\mathcal{D}}}+{\underline{h}_{\mathcal{D}}}) 
\end{array}
\end{equation} 
The function $\hat{h}_{c}$ is the central estimate of $h$, and its worst-case approximation error is the minimum among all possible estimates, considering the information that $h\in FPFS_\mathcal{D}$ \cite{TRAUB80,MILANESE2004957}. The global worst-case error pertaining to $\hat{h}_{c}$ is equal to the so-called radius of information \cite{TRAUB80}:
\begin{equation}\label{eq:worst case error central}
\max\limits_{v\in[z^{*(a1)}_{v},z^{*(a2)}_{v}]}\sup\limits_{h\in FPFS_\mathcal{D}}|h(v)-\hat{h}_{c}(v)|=\frac{1}{2}\overline{\lambda}_{\mathcal{D}}
\end{equation}
Since the worst-case approximation error obtained by $\hat{h}_{c}$ is the minimum one, this estimation approach results to be the most efficient in terms of number of samples required to guarantee the robustness condition \eqref{eq:robust requirement} for any $\mathcal{P}\in\mathcal{A}_\mathcal{D}$ (i.e., any $h\in FPFS_\mathcal{D}$).
\begin{proposition}\label{prop:central estimate}
    Consider the BOP problem  \eqref{MOP} and a tolerance vector $\boldsymbol{\delta}$. Under Assumption \ref{ass.PF}, if the central PF estimate $\hat{\mathcal{P}}_c$ \eqref{eq:central_approx} is adopted, then:    
    \begin{equation}\label{eq:main result - ROBBO}
    \begin{array}{c}
        \exists \mathcal{R}:\text{ property \eqref{eq:robust requirement} is satisfied }\forall \mathcal{P}\in FPFS_\mathcal{D}, \; \\
        \iff\\\overline{\lambda}_{\mathcal{D}}\leq2\sqrt{2}
    \end{array}
    \end{equation} 
     Moreover, the corresponding minimum number of samples is:
     \begin{equation}\label{eq.NsUni_central}
\underline{M}_c=\ceil*{\frac{1}{4}\left(\frac{\Delta_1}{\delta_1}+\frac{\Delta_2}{\delta_2}\right)}+1
\end{equation}
and the corresponding data-set $\mathcal{D}$ shall be evenly allocated in the segment $\left[z^{*(a2)}_v,z^{*(a1)}_v\right]$.
\end{proposition}
\begin{proof}
See the Appendix.
\end{proof}

Regarding the realization $\mathcal{R}$, ROBBO employs $\mathcal{R}_{SM}$, already introduced in \eqref{eq:SM realization}. The main reason is that this realization yields an additional property on the approximation error vector $\b{\varepsilon}$:
it is such that the ratio $|\varepsilon_1|/|\varepsilon_2|$ is always fixed and equal to the ratio of the prescribed error bounds: 
\begin{equation}\label{eq:balanced requirement}
    \forall \hat{\mathcal{P}},\mathcal{P}\in FPFS_\mathcal{D,\,}\forall\hat{\b{z}}\in\hat{\mathcal{P}},\,\b{z}^*=\mathcal{R}_{SM}(\hat{\b{z}}),
    \frac{|\varepsilon_1|}{|\varepsilon_2|}=\frac{\delta_1}{\delta_2} 
\end{equation}
This property can be derived along the same lines of the proof of Theorem \ref{th:2}.
When a pair $(\hat{\mathcal{P}},\mathcal{R})$ satisfies \eqref{eq:robust requirement} and \eqref{eq:balanced requirement}, we say that it is a \textit{robust and balanced} PF estimate-realization pair, see Fig. \ref{fig.setup} for visualization.

\subsection{Iterative algorithm}\label{S:Sampling}
ROBBO builds iteratively the data set $\mathcal{D}$ until the condition \eqref{eq:main result - ROBBO} is satisfied. Let us denote with $k\in\mathbb{N}$ the iteration number of the algorithm.
Besides the anchor points, computed with the auxiliary problems \eqref{eq:anchor points}, at each $k$ ROBBO computes a new sample of the PF by choosing a suitable value $\tilde{v}\in[z_v^{*(a1)},\,z_v^{*(a2)}]$ and solving the following auxiliary optimization problem:
\begin{subequations}\label{Scalarization}
\begin{gather}
    \b{x}^*=\arg\min\limits_{\b{x}\in S}{\delta_1 f_1(\b{x})+\delta_2 f_2(\b{x})}\label{Scalarization-1}\\
    \textnormal{s.t.\;}
    \delta_2f_1(\b{x})-\delta_1f_2(\b{x})-\frac{2\delta_1\delta_2}{\sqrt{2}}\tilde{v}=0\label{Scalarization-2}
\end{gather}
\end{subequations}
The corresponding PF sample is $\b{z}^*=f(\b{x}^*)$. The equality constraint \eqref{Scalarization-2} is designed to find the Pareto point $\b{z}^*$ such that its $v$-coordinate is equal to $\tilde{v}$, and its $q$-coordinate is equal to $h(\tilde{v})$. In fact, we have:
\begin{equation}\nonumber
\begin{split}
    &\begin{bmatrix}
        f_1\\f_2
    \end{bmatrix}
    =W^{-1}
    \begin{bmatrix}
        v\\q
    \end{bmatrix}
    =
    \begin{bmatrix}
       \delta_1/\sqrt{2}(v+q)\\
       \delta_2/\sqrt{2}(q-v)
    \end{bmatrix}
    \\
    & q = \frac{\sqrt{2}}{\delta_1}f_1 - v\\
    & \delta_1 f_2 = \frac{\delta_1\delta_2}{\sqrt{2}}(q-v)=\delta_2 f_1 - \frac{2\delta_1\delta_2}{\sqrt{2}}v
\end{split}
\end{equation}
Note that the gradient of the cost \eqref{Scalarization-1} is perpendicular to the gradient of the equality constraint \eqref{Scalarization-2}, both seen as functions of $f_1,\,f_2$, so that ideally the point $[f_1(\b{x}),\,f_2(\b{x})]^T$ would slide along the equality constraint until reaching the wanted PF point. Under Assumption \ref{ass.PF}, solving problem \eqref{Scalarization} always yields a PF point, both for convex and non-convex PFs. This auxiliary optimization problem also serves as a practical implementation of the realization $\b{z}^*=\mathcal{R}_{SM}(\hat{\b{z}})$ (see equation \eqref{eq:SM realization}): in this case, one shall set $\tilde{v}=\hat{z}_v$, where $\hat{z}_v$ is the $v$-coordinate of  $\hat{\b{z}}$.
We are now in position to provide a pseudo-code for ROBBO, in Algorithm 1, where the notation $\cdot^{\langle k \rangle}$ denotes a quantity that is updated at each iteration $k$. 

In virtue of Proposition \ref{prop:central estimate}, Algorithm 1 is always guaranteed to end in at most $\underline{M}_c$ iterations (including the ones to compute the anchor points), satisfying the robustness guarantees \eqref{eq:robust requirement}. In practice, a much lower number of samples may be needed, depending on the actual shape of the PF. We present some benchmarks about the actual number of iterations in Section \ref{sec:extensions, comparisons}, for different PF shapes. Regarding step e)-2) of the algorithm, if more than one value of $\tilde{v}(r)$ attain the same worst-case error, then we select the one that is farthest, according to the $v$-coordinate, from the already collected samples. If there are still more than one value with equal ranking, we pick the smallest one.

\begin{remark}
    If Assumption 1 is violated, i.e.,  the PF is discontinuous, then the auxiliary problem \eqref{Scalarization} may either be infeasible or yield a solution outside the optimal bounds, i.e. not Pareto optimal. In the first case, one can adapt the problem according to Pascoletti–Serafini scalarization to always find a Pareto optimal solution, relaxing the linear constraint \eqref{Scalarization-2}, \cite{burachikNewScalarizationTechnique2014a}. In the second, one can apply \textit{Pareto filters} to refine the data set $\mathcal{D}$, \cite{messacNormalizedNormalConstraint2003}. A rigorous extension of ROBBO to discontinuous PFs is subject of current research and outside the scope of this work.
    \end{remark}

\begin{algorithm}
\caption{ROBBO: RObust and Balanced Bi-objective Optimization}
\noindent \textbf{Inputs}: $f_1,\,f_2,\,S,\,\b{\delta}$\\
\noindent \textbf{Output}: $\hat{\mathcal{P}}_c$\\
\noindent \textbf{Algorithm}:
\begin{itemize}
    \item[a)] Set $k=1$, compute $\b{z}^{(*a1)}$ by solving \eqref{eq:anchor points-1};
    \item[b)] Set $k=2$, compute $\b{z}^{(*a2)}$ by solving \eqref{eq:anchor points-2};
    \item[c)] Initialize $\mathcal{D}^{\langle 2 \rangle}=\{\b{z}^{(*a1)},\,\b{z}^{(*a2)}\}$, compute $\Delta_1=z^{*(a2)}_1-z^{*(a1)}_1$ and $\Delta_2=z^{*(a1)}_2-z^{*(a2)}_2$, compute $\underline{M}_c$ from \eqref{eq.NsUni_central}. Define a grid of strictly increasing $v$-values $\tilde{v}^{(r)},\,r=1,\ldots,\underline{M}_c-2$ in the interval $[z_v^{*(a1)},z_v^{*(a2)}]$, equally spaced one from the next and from the anchor points' coordinates $z_v^{*(a1)},z_v^{*(a2)}$;
    \item[d)] Compute $\overline{\lambda}_{\mathcal{D}}^{\langle 2 \rangle}$ as in \eqref{eq:max global error result} with $\mathcal{D}=\mathcal{D}^{\langle 2 \rangle}$.
    \item[e)] While $\overline{\lambda}_{\mathcal{D}}^{\langle k \rangle}<2\sqrt{2}$:
    \begin{itemize}
        \item[1)] Set $k\leftarrow k+1$;
        \item[2)] Select $\tilde{v}^{\langle k \rangle}\in\arg\max\limits_{r=1,\ldots,\underline{M}_c-2}\bar{h}_{\mathcal{D}}(\tilde{v}^{(r)})-\underline{h}_{\mathcal{D}}(\tilde{v}^{(r)})$;
        \item[3)] Solve \eqref{Scalarization} with $\tilde{v}=\tilde{v}^{\langle k \rangle}$, let the obtained PF point be $\b{z}^{\langle k \rangle}$;
        \item[4)] Update (and re-order) the data set as $\mathcal{D}^{\langle k \rangle}=\mathcal{D}^{\langle k-1 \rangle}\cup \{\b{z}^{\langle k \rangle}\}$;
        \item[5)] Compute $\overline{\lambda}_{\mathcal{D}}^{\langle k \rangle}$ as in \eqref{eq:max global error result} with $\mathcal{D}=\mathcal{D}^{\langle k \rangle}$
    \end{itemize}
     \item[f)] Return $\hat{\mathcal{P}}_c$ as in \eqref{eq:central_approx} with  $\mathcal{D}=\mathcal{D}^{\langle k \rangle}$.
\end{itemize}
\end{algorithm}
\section{Extensions and comparisons}\label{sec:extensions, comparisons}

\subsection{Limited sampling budget}\label{ssec:limited budget}
A relevant situation in practice is when the DM has a limited budget of PF points to be sampled, denoted by $n_B$, and wants to know what are the attainable guaranteed worst-case error bounds on each of the two objective functions, $\bar{\delta}_1,\,\bar{\delta}_2$, that comply with a wanted ratio $\alpha=\frac{\bar{\delta}_1}{\bar{\delta}_2}$, representing the relative importance of the two criteria. Using the central estimate (which guarantees the smallest worst-case error) and evenly distributed samples as done by ROBBO (which guarantee satisfaction of property \eqref{eq:robust requirement} with the smallest number of samples), we can impose the following condition:
\begin{equation}\nonumber
    2\overline{\lambda}_{\mathcal{D}}=\frac{V_a}{n_B-1}=\frac{1}{n_B-1}\frac{\sqrt{2}}{2}\left(\frac{\Delta_1}{\bar{\delta}_1}+\frac{\Delta_2}{\bar{\delta}_2}\right)
\end{equation}
which stems from Proposition \ref{prop:central estimate} and includes the computation of the two anchor points as part of the budget $n_B$. Then, using $\bar{\delta}_1=\alpha\bar{\delta}_2$, we obtain:
\begin{equation}\label{eq:inverse_prob-2}
\bar{\delta}_1 = \frac{\Delta_1+\alpha\Delta_2}{4(n_B-1)},\;\;
\bar{\delta}_2 = \frac{\Delta_1+\alpha\Delta_2}{4\alpha(n_B-1)}
\end{equation}
This result can be used to pre-compute, after obtaining the anchor points, the attainable worst-case accuracy guarantees as a function of the sampling budget $n_B$ to be allocated. A case study for this approach is presented in Section \ref{Sec:Stirred-Tank Reactor}.

\subsection{Linear interpolation}\label{sec:linear interp}
\begin{figure}[hbt!]
\centerline{\includegraphics[trim=150 185 419 57,clip,scale=0.65]{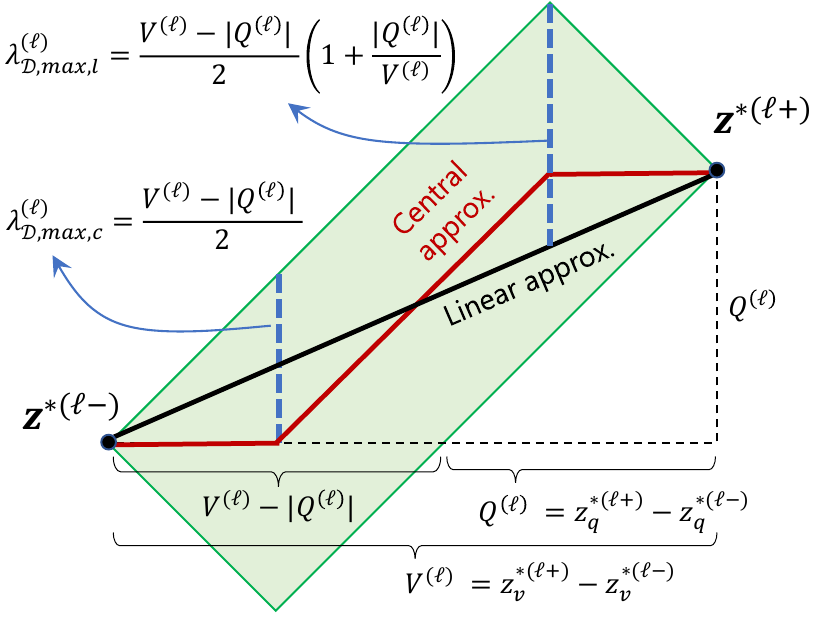}}
\caption{Computation of $\lambda_{\mathcal{D},max,c}^{(\ell)},\,\lambda_{\mathcal{D},max,l}^{(\ell)}$, for the central and linear approximation, in the generic interval $R^{(\ell)}$}
\label{fig.MaxUn}
\end{figure}

An alternative estimation technique, widely used in practice, is the piecewise linear approximation, denoted by $\hat{h}_l(v)$. For a generic interval $R^{(\ell)}=[\b{z}^{*(\ell-)},\,\b{z}^{*(\ell+)}]\subset [z_v^{*(a1)},\,z_v^{*(a2)}]$
(see Section \ref{SS:minimum M}), this estimate is given by
\begin{equation*}
    \hat{h}_l(v)=z^{*(\ell-)}_{q}+\frac{z^{*(\ell+)}_{q}-z^{*(\ell-)}_{q}}{z^{*(\ell+)}_{v}-z^{*(\ell-)}_{v}}(v-z^{*(\ell-)}_{v}),\,\forall v\in R^{(\ell)}
\end{equation*}
For simplicity, we adopt a graphical interpretation to study the guaranteed accuracy properties of this estimation technique and compare them to those of the central estimate, see Fig. \ref{fig.MaxUn}. 
We denote with $\lambda_{\mathcal{D},max,c}^{(\ell)}$ and with $\lambda_{\mathcal{D},max,l}^{(\ell)}$ the worst-case error pertaining to the central and linear approximation, respectively. Regarding the central estimate, consistently with Proposition \ref{prop: max global error} and \eqref{eq:worst case error central} we have
$\lambda_{\mathcal{D},max,c}^{(\ell)}=\frac{V^{(\ell)}-|Q^{(\ell)}|}{2}$
and this value is attained for all $v$-values in the interval $[v'',\,v']$ (see Fig. \ref{fig.MaxUn} and the proof of Proposition \ref{prop: max global error}). About the linear estimate, geometric considerations lead to 
$\lambda_{\mathcal{D},max,l}^{(\ell)}=\frac{V^{(\ell)}-|Q^{(\ell)}|}{2}\left(1+\frac{|Q^{(\ell)}|}{V^{(\ell)}}\right)$. 
It is immediate to note that the two worst-case errors are the same only if the $q-$coordinates of $\b{z}^{*(\ell-)},\,\b{z}^{*(\ell+)}$ coincide, i.e., $|Q^{(\ell)}|=0$. This is the situation when the worst-case error is largest in the interval $R^{(\ell)}$, as argued in the proof of Proposition \ref{prop:minimum M}. Otherwise, the central estimate has generally better worst-case performance, depending on the difference $V^{(\ell)}-|Q^{(\ell)}|$. Note that, due to the Lipschitz continuity property of the rotated PF with constant smaller than one (Remark  \ref{rem:Lipschitz}), we always have  $V^{(\ell)}>|Q^{(\ell)}|$. Finally, note that, for the linear approximation, the worst-case error is attained only at coordinates $v''$ and $v'$ and not in-between them.
\subsection{Greedy bisection sampling strategy}\label{sec:max_un}
It is interesting to compare the sampling approach adopted in ROBBO to a more common,  greedy bisection method. In the latter, instead of picking a point from the uniform grid, at each iteration one selects, as value of $\tilde{v}$ in problem \eqref{Scalarization}, the midpoint with the largest uncertainty among all intervals $R^{(\ell)}$, thus bisecting that interval. From the considerations of Section \ref{sec:linear interp}, the midpoint always corresponds to the largest worst-case error $\lambda_{\mathcal{D},max,c}^{(\ell)}$ for the central estimate, thus sampling it can generally bring a significant reduction of the uncertainty. For the linear approximation, to pick the point with maximum uncertainty one could sample at coordinate $v'$ or $v''$, i.e. $ z^{*(\ell)+}_v-(V^{(\ell)}-|Q^{(\ell)}|)/2$ or $z^{*(\ell-)}_v+(V^{(\ell)}-|Q^{(\ell)}|)/2$.\\
If the greedy bisection strategy is adopted, the algorithm still converges in finite iterations and we can pre-compute their maximum number. This corresponds to $\underline{M}_c$ only in a very specific case, otherwise it is larger. In fact, if the bisection approach is used, in the worst case (i.e. when the actual PF is such that its $q$-coordinates are all equal, so that $|Q^{(\ell)}|=0,\,\forall \ell$),
the global worst-case error bound $\max\limits_{\ell}\lambda_{\mathcal{D},max,c}^{(\ell)}$ is halved only when all intervals have been halved. Each time that this happens, 
we say that one epoch has passed. We can then compute the worst-case number of epochs $n_e$ and the upper bound of samples, $\underline{M}_g$, for such a greedy strategy
to guarantee satisfaction of \eqref{eq:main result - ROBBO}, as a function of the distance $V_a$ between the anchor points' $v$-coordinates:
\begin{equation}\nonumber
\begin{split}
      &\frac{V_a}{2^{n_e}}\leq2\sqrt{2}\Rightarrow n_e = \ceil*{\log_{2}\left(\frac{V_a}{2\sqrt{2}}\right)}\\
      &\underline{M}_g = 2+\sum^{n_e}_{i=1}{2^{(i-1)}}
\end{split}
\end{equation}
where $\underline{M}_g$ takes into account the anchor points. We note that $\underline{M}_g$ of the greedy bisection strategy is equal to  $\underline{M}_c$ of the uniform sampling only in the very special case that the ratio $V_a/(2\sqrt{2})$ is a power of $2$, otherwise, it is larger.

\subsection{Empirical benchmark of estimation and sampling strategies}
\begin{figure}
\centerline{\includegraphics[trim=312 5 200 5,clip,scale=0.55]{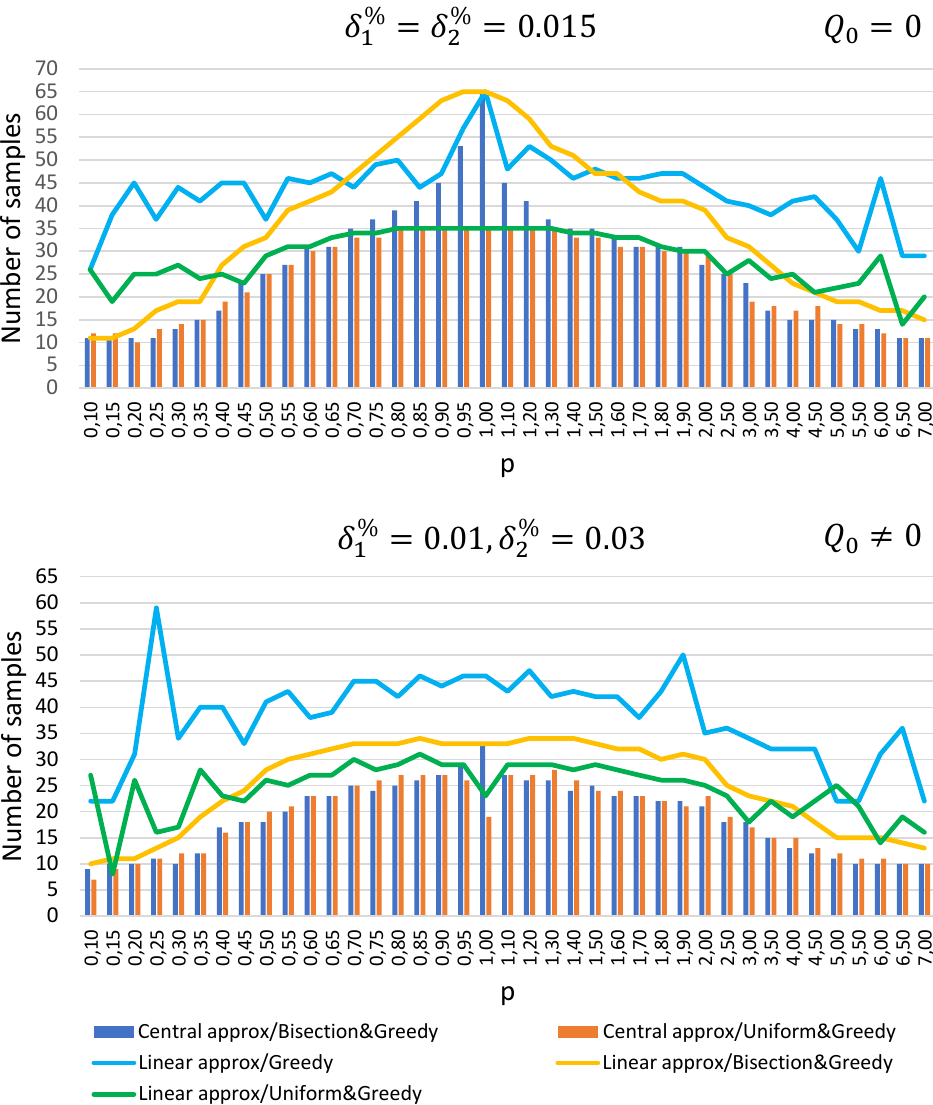}}
\caption{Empirical comparison of the iterative sampling heuristics when varying the front shape}
\label{fig.heuristics}
\end{figure}

The theoretical findings of Sections \ref{sec:linear interp} and \ref{sec:max_un} refer to worst-case situations. It is of interest to evaluate, via benchmarks that are representative enough of different possible cases, what are the actual performance of using either the central or the linear estimation and either a uniform distribution of samples or the greedy bisection approach. Here, by performance we refer to the actual number of samples required to eventually satisfy the robust accuracy requirement \eqref{eq:robust requirement} with the chosen estimation method.\\
We perform the analysis comparing the following alternatives:
\begin{itemize}
    \item central estimate with uniformly allocated points (i.e., ROBBO);
    \item central estimate with greedy bisection sampling;
    \item linear interpolation with uniform distribution of points;
    \item linear interpolation with greedy bisection sampling;
    \item linear interpolation with greedy sampling of the points with largest uncertainty (i.e. points with coordinates $v'$ or $v''$ mentioned in Section \ref{sec:max_un}).
\end{itemize}
We set the elements of the tolerance vector as percentages of the ranges $\Delta_1,\,\Delta_2$, i.e. $\delta_1=\delta_1^{\%}\Delta_1,\delta_2=\delta_2^{\%}\Delta_2$, and consider the following BOP:
\begin{subequations}\label{eq:BOPtest}
\begin{gather}
        \min_{\b{x}}\b{f}(\b{x})=[x_1,x_2]^T\\
        \text{s.t.}\nonumber\\
        {\left(\sum_{i=0}^2{x_i}^{p}\right)}^{\frac{1}{p}}\geq10\label{eq:BOPtest_const}\\
        x_1,x_2\geq0
\end{gather}
\end{subequations}
The PF of this problem is given by $f_2=(10^p-{f_1}^p)^{\frac{1}{p}}$, derived from constraint (\ref{eq:BOPtest_const}). The anchor points are $\b{z}^{*(a1)}=[0, 10]^T$, $\b{z}^{*(a2)}=[10, 0]^T$. The PF is convex for $0<p<1$, it is linear for $p=1$, and concave for $p>1$. In addition to the shape, the value of $Q_0=|z^{*(a2)}_q-z^{*(a1)}_q|$ can be controlled by varying $\delta_1^{\%},\,\delta_2^{\%}$: since $\Delta_1=z^{*(a2)}_1,\Delta_2=z^{*(a1)}_2$, we have in fact $Q_0=\frac{\sqrt{2}}{2}\left(\frac{1}{\delta_1^{\%}}-\frac{1}{\delta_2^{\%}}\right)$. Thus, when $\delta_1^{\%}=\delta_2^{\%}$ the two anchor points have the same $q$-coordinate. This allows us to also test the worst-case situation, by setting $p=1$ (linear front) and $\delta_1^{\%}=\delta_2^{\%}$ (the whole rotated front with the same $q$-coordinate).\\
We present the results obtained in case of $Q_0=0$ with $\delta_1^{\%}=\delta_2^{\%}=1.5\%$, and of $Q\neq0$ with $\delta_1^{\%}=1\%,\;\delta_2^{\%}=3\%$. In both cases, we spanned different $p\in[0.01,7]$. The results are reported in Fig. \ref{fig.heuristics}.
\noindent Overall, we see that the central estimate performs better than the linear interpolation and, as expected from the theory, they behave similarly only with uniformly spaced samples when the shape of the front is close to linear (i.e., when $p\approx 1$ in the upper plot).
This comparison confirms that the central estimate combined with uniformly spaced samples, as implemented in ROBBO, is overall the best option.
A demo implementation of this example is openly available at \cite{ROBBOweb}.

\subsection{Theoretical comparison with popular PF estimation methods}
\begin{figure}
\centerline{\includegraphics[trim=80 245 357.5 75,clip,scale=0.47]{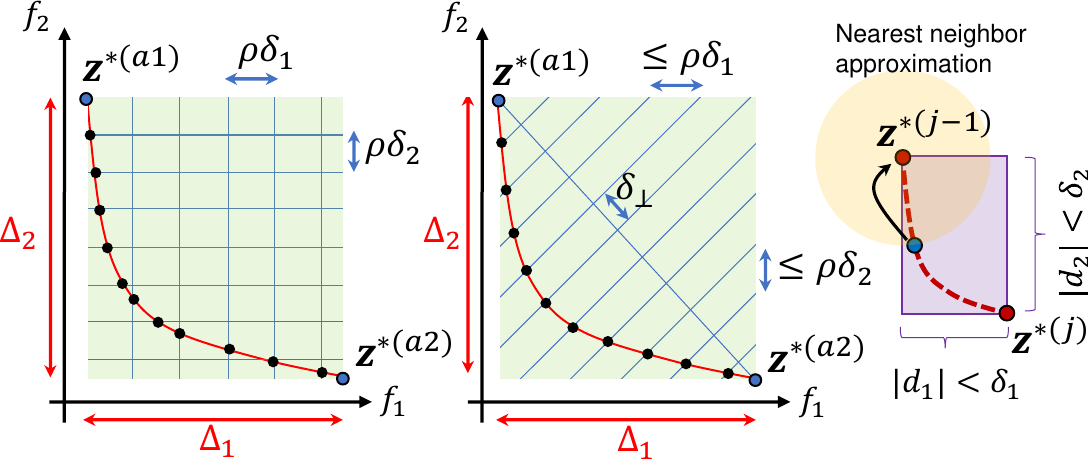}}
\caption{Uniform sampling for EC (left plot), NBI (central plot) to satisfy the requirement \eqref{eq:robust requirement}, and nearest neighbor PF approximation (right)}
\label{fig.EC_NBI}
\end{figure}

We now compare the worst-case required number of samples obtained by $ROBBO$ with that of some of the most popular exact strategies for PF approximation. For this purpose, we have selected the $\epsilon$-constraint approach \cite{CHIANDUSSI2012912}, denoted by EC, and the Normal Boundary Intersection one \cite{dasNormalBoundaryIntersectionNew1998}, NBI. 
Because these methods deliver a point-wise approximation of the PF (i.e., nearest neighbor), we consider the distance between a non-computed point and its nearest neighbor as the approximation error (see Fig. \ref{fig.EC_NBI}, right). Thus, to satisfy the robust accuracy condition \eqref{eq:robust requirement}, the distance between two sampled points $\b{d}=\b{z}^{*(j+1)}-\b{z}^{*(j)}$ must satisfy:
\begin{equation}\label{eq:conditionNN}
    |d_1|\leq\delta_1\wedge|d_2|\leq\delta_2
\end{equation}
Using the EC method, a point-wise approximation of PF is provided by solving several instances of the scalar problem that minimizes one of the two objectives, while constraining the other objective to be less than or equal to a value $\epsilon$ \cite{CHIANDUSSI2012912}. At each new instance, $\epsilon$ is increased by a fixed quantity to obtain a uniform distribution of points with respect to the constrained criterion. To satisfy condition \eqref{eq:conditionNN}, one needs to repeat the procedure for both objectives
(see Fig. \ref{fig.EC_NBI}, left).
The corresponding worst-case number $\underline{M}_{EC}$ of samples is therefore:
\begin{equation}\label{eq:Ns_EC}
    \begin{split}
        \underline{M}_{EC} = \ceil*{\frac{\Delta_1}{\delta1} + \frac{\Delta_2}{\delta2}} +1
    \end{split}
\end{equation}

\noindent The NBI method \cite{dasNormalBoundaryIntersectionNew1998} produces samples that are uniformly spaced when projected on the segment connecting the anchor points, see Fig. \ref{fig.EC_NBI} (center). We denote with $W$ such a segment, and with
$\delta_{\perp}$ the distance between the projections of two subsequent samples on $W$.
Let us consider the quantities $l_1=\delta_1\cos{(\beta)}$ and $l_2=\delta_2\sin{(\beta)}$, where $\beta =\tan^{-1}(\Delta_2/\Delta_1)$. In the worst case, to fulfill \eqref{eq:conditionNN} it must hold 
$\delta_{\perp}\leq\min{(l_1,l_2)}$. 
Then, we can compute the worst-case number of samples $\underline{M}_{NBI}$ by dividing the length of $W$ by $\delta_{\perp}$ that satisfies this condition.
Without loss of generality, we consider $\delta_1 = \alpha\delta_2$ and $\Delta_1=\gamma\Delta_2$. Then, assuming $l_2<l_1$ we have:

    \begin{align*}
        &\frac{\Delta_2}{\delta_{\perp}\sin{(\beta)}} = \frac{\Delta_2}{\delta_2 {(\sin{(\beta)})}^2}= \frac{\Delta_2}{\delta_2}\frac{1+{\tan(\beta)}^2}{{\tan(\beta)}^2}=\\
        &= \frac{\Delta_2}{\delta_2}\frac{1+{(\Delta_2/\Delta_1)}^2}{{(\Delta_2/\Delta_1)}^2} = \frac{{\Delta_1}^2+{\Delta_2}^2}{\delta_2\Delta_2}=\\
        &=\left(\frac{{\Delta_1}^2}{\delta_2\Delta_2}+\frac{\Delta_2}{\delta_2}\right)=\left(\frac{{\Delta_1}^2}{\frac{\delta_1}{\alpha}\frac{\Delta_1}{\gamma}}+\frac{\Delta_2}{\delta_2}\right) =\left(\frac{{\Delta_1}}{\delta_1}\alpha\gamma+\frac{\Delta_2}{\delta_2}\right) \\
         &l_2<l_1\Rightarrow\delta_2\sin{(\beta)}<\delta_1\cos{(\beta)}\Rightarrow\delta_2\sin{(\beta)}<\alpha\delta_2\cos{(\beta)}\\
         &\Rightarrow \alpha>\frac{\sin{(\beta)}}\Rightarrow{\cos{(\beta)}}\Rightarrow \alpha>\frac{\Delta_2}{\gamma\Delta_2}\Rightarrow\alpha\gamma>1
    \end{align*}
from which we derive: 
\begin{equation}\label{eq:Ns_NBI1}
    \underline{M}_{NBI} = \ceil*{\frac{{\Delta_1}}{\delta_1}\alpha\gamma+\frac{\Delta_2}{\delta_2}} +1, \alpha\gamma>1
\end{equation}
Assuming $l_1<l_2$, with similar passages we obtain:
\begin{equation}\label{eq:Ns_NBI2}
    \underline{M}_{NBI} = \ceil*{\frac{{\Delta_1}}{\delta_1}+\frac{\Delta_2}{\delta_2}\frac{1}{\alpha\gamma}} +1, \frac{1}{\alpha\gamma}>1
\end{equation}
Finally if  $l_1=l_2$ we find:
\begin{equation}\label{eq:Ns_NBI3}
    \underline{M}_{NBI} = \underline{M}_{EC} = \ceil*{\frac{{\Delta_1}}{\delta_1}+\frac{\Delta_2}{\delta_2}} +1
\end{equation}
Now, comparing $\underline{M}_c$ \eqref{eq.NsUni_central}, $\underline{M}_{EC}$ \eqref{eq:Ns_EC}, and $\underline{M}_{NBI}$ \eqref{eq:Ns_NBI1}-\eqref{eq:Ns_NBI3}, we see that
NBI needs, in the worst case, at least
the same number of samples as EC, which is approximately four times larger than in our method, where we exploit the optimal bounds and a continuous estimated PF instead of a nearest-neighbor one.
\section{Case Studies}\label{S:Results}
\subsection{Path-following positioning system}
\begin{figure}
\centerline{\includegraphics[trim=178 115 210 130,clip,scale=0.55]{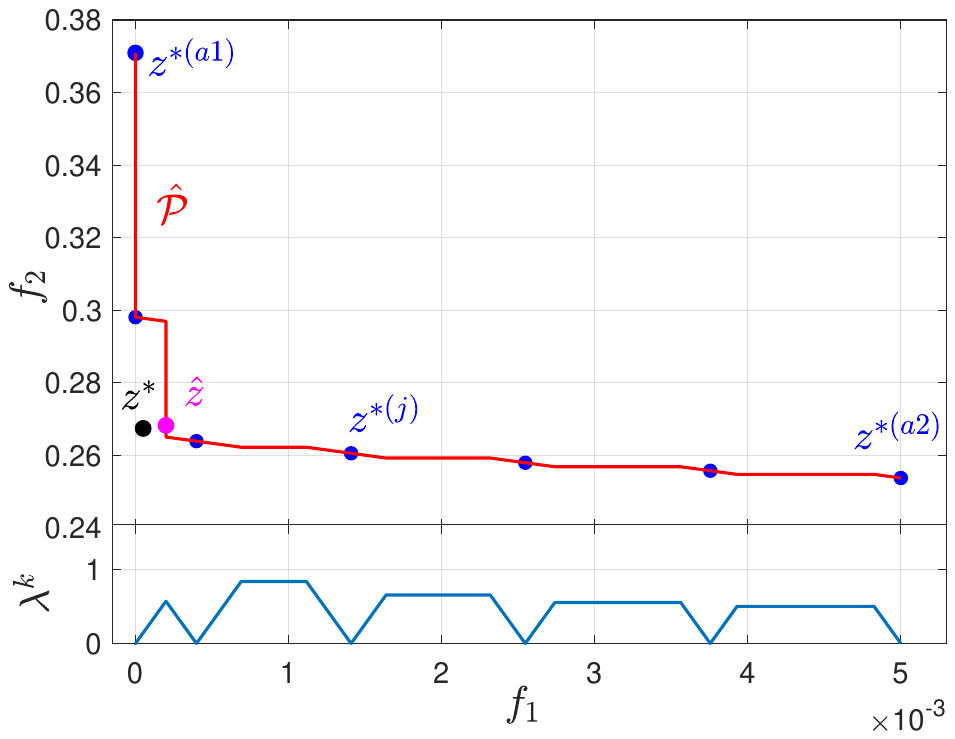}}
\caption{Upper plot: central estimate of the Pareto front. Lower plot: uncertainty relative to the approximation, represented along the $f_1$ axis}
\label{fig.PFresult}
\end{figure}

\begin{figure}
\centerline{\includegraphics[trim=115 90 130 125,clip,scale=0.425]{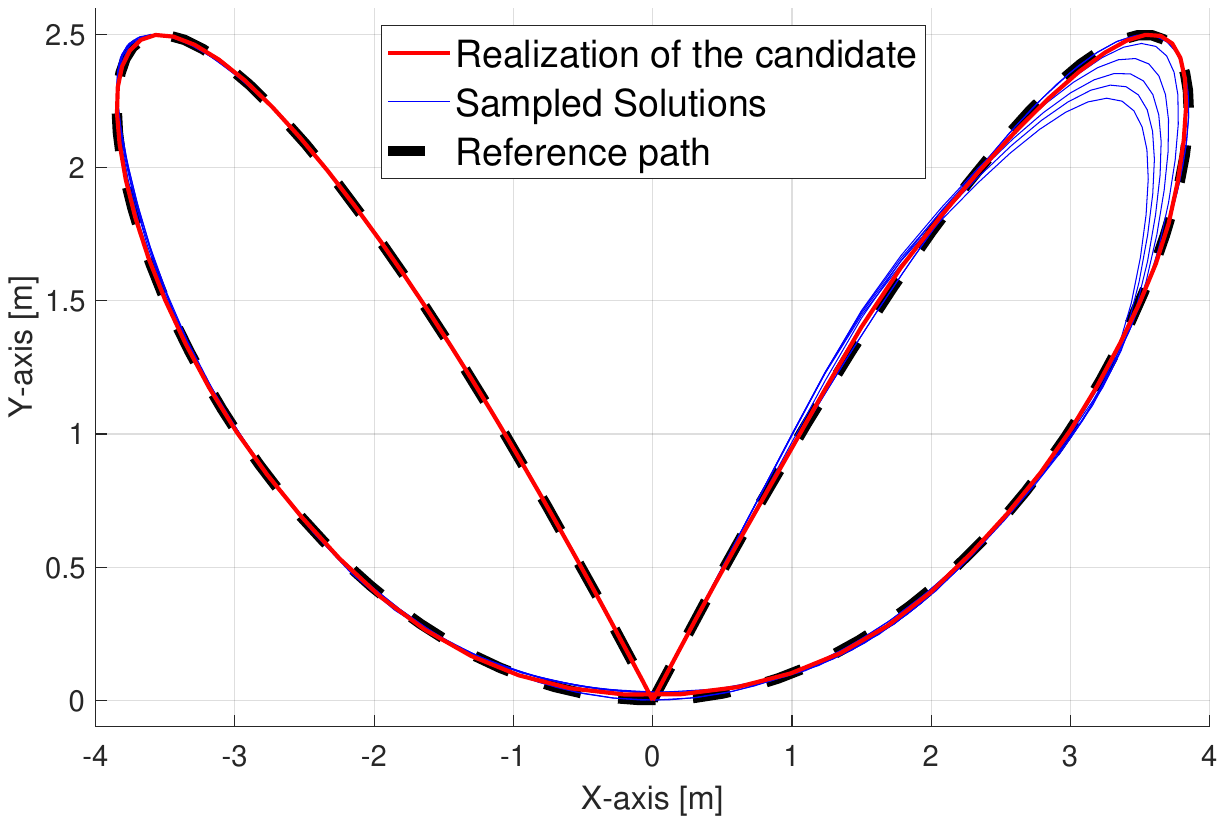}}
\caption{Path following: reference path and all Pareto optimal solutions computed by ROBBO}
\label{fig.trajectory}
\end{figure}

We present the application of ROBBO to a constrained-path following problem for a 2-axis positioning system, formulated as a Finite Horizon Optimal Control Problem (FHOCP). The system model in continuous time $\tau$ is:
\begin{equation*}
\dot{\boldsymbol{\xi}}(\tau)=
\left[
\begin{smallmatrix}
	0&1&0&0\\
	0&-\beta_f/m&0&0\\
	0&0&0&1\\
	0&0&0&-\beta_f/m
\end{smallmatrix}
\right]
\boldsymbol{\xi}(\tau)+
\left[
\begin{smallmatrix}
	0&0\\1/m&0\\0&0\\0&1/m\\
\end{smallmatrix}
\right]
\b{u}(\tau)
\end{equation*}
The states $\xi_1$ and $\xi_3$ are the horizontal and vertical position of the end effector (e.e.) in m, $\xi_2$ and $\xi_4$ are the e.e. horizontal and vertical speed in m/s, $\beta_f=0.1\,$Ns/m is the viscous friction coefficient, $m=2\,$kg is the mass of the moving parts, $u_1$ and $u_2$ are the horizontal and vertical actuator forces in N.
The reference path is denoted by $\bar{\b{\xi}}(\theta)$, where  $ \theta\in[0,1]$, is a curvilinear coordinate such that $\theta(0)=0$ indicates the starting point and $\theta(H)=1$ the completion of the path.\\
Adopting a time discretization, we express the optimal control problem as a BOP, where the two conflicting objectives are related to the tracking accuracy and the advancement of the reference coordinates:
\begin{align*}
&f_1(\boldsymbol{\xi}(t),\theta(t))=\frac{1}{H}\sum\limits_{t=0}^{H}(\bar{\boldsymbol{\xi}}(\theta(t))-\boldsymbol{\xi}(t))^TL(\bar{\boldsymbol{\xi}}(\theta(t))-\boldsymbol{\xi}(t))\\
&f_2(\boldsymbol{\xi}(t),\theta(t))=\frac{1}{H}\sum\limits_{t=0}^{H}(1-\theta(t))^2
\end{align*}
where $t=1\dots,H$ is the discrete time index, $H$ the terminal time defined by the user, and $L=\text{diag}([1,0,1,0])$.
We formulate the following FHOCP:
\begin{subequations}\label{eq:cnoec-specific-path-following}
\begin{gather}
\min\limits_{U,\Delta\Theta}\,f_1(\boldsymbol{\xi}(t),\theta(t)),\,f_2(\boldsymbol{\xi}(t),\theta(t))\label{eq:cnoec-specific-path-following-cost}\\
\text{s.t.}\nonumber\\
\boldsymbol{\xi}(t+1)=f_{\boldsymbol{\xi}}(\boldsymbol{\xi}(t),\b{u}(t)),\;t=0,\ldots,H-1\label{eq:cnoec-specific-path-following-dynamics}\\
\theta(t+1)=\theta(t)+\Delta\theta(t),\;t=0,\ldots,H-1\label{eq:cnoec-specific-path-following-time}\\
\Delta\theta(t)\geq0,\;t=0,\ldots,H-1\label{eq:cnoec-specific-path-following-constraints-theta}\\
-\overline{\b{u}}\geq\b{u}(t)\geq\overline{\b{u}},\;t=0,\ldots,H-1\label{eq:cnoec-specific-path-following-constraints}\\
\boldsymbol{\xi}(0)=[0,0,0,0]^T\wedge\theta(0)=0\wedge\theta(H)=1\label{eq:cnoec-specific-path-following-terminal time}
\end{gather}
\end{subequations}
where $\Delta\Theta=[\Delta\theta(0),\ldots,\Delta\theta(H-1)]^T$ and
$\b{U}=[\b{u}(1)^T,\ldots\allowbreak \ldots,\b{u}(H-1)^T]^T$. Constraint \eqref{eq:cnoec-specific-path-following-dynamics} enforces the system dynamics, according to a discretized model of the continuous LTI system
with the trapezoid method with sampling period
$T_s=0.05s$.  Constraints \eqref{eq:cnoec-specific-path-following-time} and \eqref{eq:cnoec-specific-path-following-constraints-theta} account for the progress status along the path,
while \eqref{eq:cnoec-specific-path-following-constraints} defines force limits, with $\overline{\boldsymbol{u}}=[20, 20]^T\,$N, finally \eqref{eq:cnoec-specific-path-following-terminal time} sets the initial state and initial and terminal condition on the path progress. We set the control horizon $H=120$ steps, corresponding to a time window of 6 seconds.\\
The reference path is defined by polynomials in $\theta$:
$\xi_1(\theta)=
80\theta^3-120\theta^2+40\theta$, $\xi_2(\theta)=
-160\theta^4+320\theta^3-200\theta^2+40\theta$. 
Moreover, we introduce an additional constraint $f_1\leq\num{5e-3}$, to restrict the PF estimate to the most interesting portion, excluding solutions with trajectories too far from the reference. The two computed anchor points are $\b{z}^{*(a1)}=[10^{-8}, 3.71\,10^{-1}]$ and $\b{z}^{*(a2)}=[5\,10^{-3}, 2.54\,10^{-1}]$, and the objectives' tolerances are set as $\delta_1 = 3.5\,10^{-4}$ and $\delta_2 =  2\,10^{-3}$ (i.e. approx. $7\%$ and $1.7\%$ of the range of values given by the anchor points).\\
From the anchor points computation we derive $\underline{M}_c=20$, by (\ref{eq.NsUni_central}), and ROBBO reaches the accuracy guarantees in 5 iterations after computing the anchor points, i.e. $k=7$ samples in total. Fig. \ref{fig.PFresult} shows that the algorithm sampled fewer points where $f_1$ is almost zero (i.e., very small tracking error) with large completion time, since the corresponding solutions are practically equivalent from the perspective of $f_1$.
Then, taking the role of the DM, we selected a candidate solution on the approximate PF, close to the knee, i.e. the point at which the completion time is greatly reduced but with a very low-cost w.r.t. tracking error.  The selected candidate is $\hat{\b{z}}= [1.99\,10^{-4}, 2.683\,10^{-1}]$, and its corresponding realization $\b{z}^*$ is $\mathcal{R}_{SM}(\hat{\b{z}})=[5\,10^{-5}, 2.674\,10^{-1}]$, both displayed in Fig. \ref{fig.PFresult}. The figure also shows the value of the error bound along the approximation, in the graph below the $\hat{\mathcal{P}}$ plot.
The simulation of the plant behavior corresponding to the chosen PF point is provided in Fig. \ref{fig.trajectory}, where the dashed line is the path to be tracked, the red line is the simulation of the realization and the blue lines are the trajectories related to the PF points sampled by ROBBO to estimate the front. Comparing the solution chosen by the DM and its realization, according to the theory we obtained an approximation error $\b{\varepsilon}=[1.49\,10^{-4}, 8.5\,10^{-4}]$, such that $\varepsilon_1<\delta_1\land\varepsilon_2<\delta_2$ and $\frac{\varepsilon_1}{\varepsilon_2}=\frac{\delta_1}{\delta_2}=0.175$. 
\subsection{Continuous Stirred-Tank Reactor}\label{Sec:Stirred-Tank Reactor}
We consider the model of a Continuous-flow Stirred Tank Reactor (CSTR) with an exothermic, irreversible reaction $A\rightarrow B$ (see e.g. \cite{MAGNI20011351}):
\begin{equation*}
    \begin{array}{rcl}
        \dot{C}_A(\tau) &=  &\dfrac{F(C_0-C_A(\tau))}{10^3 \pi r^2 h_t}-k_0C_A(\tau)e^{-\dfrac{E}{R T(\tau)}}\\
         \dot{C}_B(\tau) &=  &\dfrac{-F(C_B(\tau))}{10^3 \pi r^2 h_t}+k_0C_A(\tau)e^{-\dfrac{E}{R T(\tau)}}\\
         \dot{T}(\tau) &=  &\dfrac{F(T_0-T(\tau))}{10^3 \pi r^2 h_t}-\dfrac{\Delta H}{\rho C_p}k_0C_A(\tau)e^{-\dfrac{E}{R T(\tau)}}\\
         &&+\dfrac{2U_e}{r\rho C_p}(T_c(\tau)-T(\tau))
    \end{array}
\end{equation*}
\begin{figure}
\centerline{\includegraphics[trim=60 235 220 237,clip,scale=0.8]{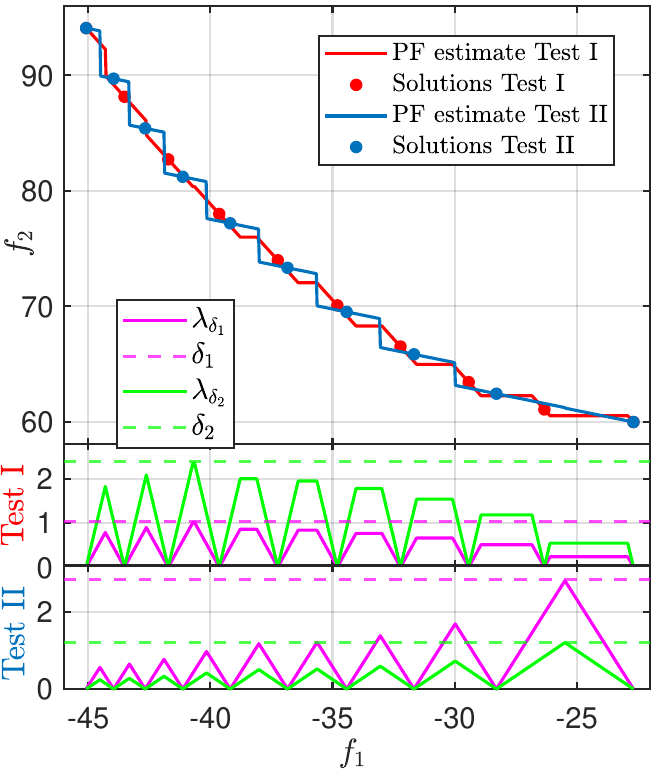}}
\caption{CSTR: Test I and test II results. Estimated PF of the economically optimal steady states for a Continuous-flow Stirred Tank Reactor (upper plot), worst-case approximation error (lower plots)}
\label{fig.3&7}
\end{figure}
\begin{figure}
\centerline{\includegraphics[trim=75 255 90 278,clip,scale=0.59]{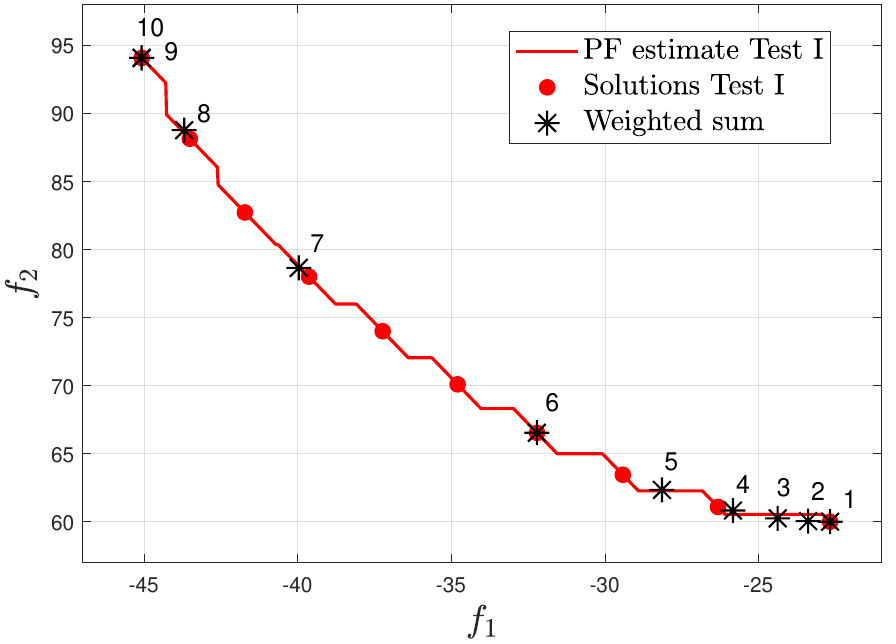}}
\caption{CSTR: comparison between Test I and scalarization by convex combination of the two objectives, with the same budget $n_B=10$ of samples.}
\label{fig.compare scalarization}
\end{figure}
Where $C_A,\,C_B,\,T$ are the concentrations of reactant $A$, product $B$, and the mixture temperature; $T_c$ the cooling temperature, $F$ the flow rate, $C_0$  the concentration of $A$ in the input flow, $r,\,h_t$ the reactor's radius and level, $k_0,\,E$ the reaction frequency coefficient and activation energy, $R,\,\Delta H,\,U_e,\,C_p$ the gas constant, reaction enthalpy, reactor' surface heat conductance, and heat capacity. We used the parameters $F=100\,$l/min, $r=0.219\,$m, $k_0=7.2 10^{10}\,$min$^{-1}$, $R=8.314\,$J/(mol K), $E=7.27475 10^{4}$J/mol, $U_e=5.49 10^{4}\,$J/(min m$^2$ K), $\rho=1\,$ kg/l, $C_p=2.4 10^2\,$J/(kg K), $\Delta H=-5\,10{^4}\,$J/(mol l),   $T_0=350\,$K, $C_0=1\,$mol/l.\\
\noindent We employ ROBBO to estimate the Pareto Front of steady-states for this system, i.e. values of $\boldsymbol{x}=[C_A,\,C_B,\,$ $T,\,T_c,\,F,\,h_t]^T$ such that $\dot{C}_A$, $\dot{C}_B$, $\dot{T}$=0, subject to the operational constraints $C_{A,B}\in[0.1,\,1]\,$mol/l, $T\in[290,\,340]\,$K, $T_c\in[280,\,350]\,$K, $F\in[80,\,120]\,$l/min, $h_t\in[0.1,\,1]\,$m  and considering these two objectives:
\begin{equation*}
    \begin{array}{rcl}
         f_1(x)& = &-C_B F\\
         f_2(x)& = &0.75 F +0.01(T_c-300)^2
    \end{array}
\end{equation*}
$f_1$ is the production rate of $B$, while $f_2$ is the total production cost, i.e. cost of reactants plus cost for the thermal regulation via the cooling jacket.\\
We provide ROBBO  with a budget of $n_B=10$ samples including the anchor points, see Section \ref{ssec:limited budget}, and perform two tests. In test I, we set the tolerance ratio $\alpha= \frac{\bar{\delta}_1}{\bar{\delta}_2}=3/7$ to favor better accuracy for $f_1$, while in test II we set $\alpha=7/3$.
When adopting this strategy, step c) of algorithm 1 presents the following changes: $\underline{M}_c$ is not computed but directly set equal to $n_B$, $\delta_1=\bar{\delta}_1, \delta_2=\bar{\delta}_2$ are computed using \eqref{eq:inverse_prob-2}. Then, at step e) the algorithm is iterated as long as $k\leq n_B$, thus exploiting the whole sampling budget.
The results are displayed in Fig. \ref{fig.3&7}, where the upper plot presents the PF estimates for both test I and test II, in red and blue, respectively.
The two lower plots present the worst-case realization errors for each objective individually, computed as $\lambda_{\delta_1}\left(\hat{z}_1^{(j)}\right)=\bar\delta_1\frac{\sqrt{2}}{2}\lambda_D\left(\hat{z}_v^{(j)}\right)$ and $\lambda_{\delta_2}\left(\hat{z}_1^{(j)}\right)=\bar\delta_2\frac{\sqrt{2}}{2}\lambda_D\left(\hat{z}_v^{(j)}\right)$.
In test I we obtain the guaranteed tolerances $\bar{\delta}_1=1.0284, \bar{\delta}_2=2.3997$, while for test II we obtain $\bar{\delta}_1=2.8314, \bar{\delta}_2=1.2134$. As expected from the theory, their ratios are equal to the imposed $\alpha$ values. Furthermore, as shown in Fig. \ref{fig.3&7}, the guaranteed error bounds $\lambda_{\delta_1}$ and $\lambda_{\delta_2}$ are always below the precomputed values of $\bar{\delta}_1$ and $\bar{\delta}_2$. It is interesting to notice how the samples' positions in terms of  $(f_1,f_2)$-values change between the two tests, while in the $(v,q)$-coordinates the approach always adopts an even distribution of samples. This is due to the effect of the tolerance values on the linear transformation (\ref{eq:transformation}), recalling that the identification is performed in coordinates $(v,q)$.
Finally, Fig. \ref{fig.compare scalarization} presents a comparison between Test I and the samples obtained by minimizing the convex combination $\beta f_1+(1-\beta)f_2$ with the same budget $n_B=10$ of $\beta$ values evenly distributed in the interval $[0,1]$. Such a scalarization method is a rather commonly-used approach to attempt the PF approximation. It can be clearly noted that the scalarization results are extremely poor with respect to ROBBO, with a large part of the PF without any samples and corresponding large approximation error.

\section{Conclusions and future developments}\label{S:conclusions}
We presented new results in the field of bi-objective optimization, delivering the algorithm ROBBO, which builds an estimate of the Pareto Front with guaranteed accuracy, low number of required samples, and a wanted balancing between the estimation errors with respect to the two criteria. We proved these properties theoretically and compared the approach with other ones, both theoretically and empirically, showing the superiority of ROBBO in terms of the required number of samples to attain the wanted guarantees.
Future research will focus on problems with more than two objectives, discontinuous fronts, and uncertain fronts, i.e. with objective function outcomes affected by some intrinsic variability.

\section*{Appendix}
\textit{Proof of Theorem \ref{T:bounds}}. Consider a value $v\in[z_v^{*(a1)},z_v^{*(a2)}]$, and take any $z^{*(j)}\in\mathcal{D}$. Assume that $z_v^{*(j)}<v$, i.e. $v=z_v^{*(j)}+\Delta v$ for some $\Delta v>0$. Now consider a value $q=z_q^{*(j)}+|v-z_v^{*(j)}|+\Delta q$, where $\Delta q\in\mathbb{R}$. Since $\Delta v>0$, we have that $q=z_q^{*(j)}+\Delta v+\Delta q$. Then, consider the point $[v,q]^T$ and apply the inverse transformation:
\[
\b{z}=W^{-1}
\left[
\begin{array}{c}
v\\q
\end{array}\right]=\frac{\delta_1\delta_2}{\sqrt{2}}\left[
\begin{array}{r}
\frac{v}{\delta_2}+\frac{q}{\delta_2}\\-\frac{v}{\delta_1}+\frac{q}{\delta_1}
\end{array}\right]
\]
By substituting $v=z_v^{*(j)}+\Delta v$ and $q=z_q^{*(j)}+\Delta v+\Delta q$ and after a few manipulations, we obtain:
\[
\b{z}=
\b{z}^{*(j)}+\left[
\begin{array}{l}
\frac{\delta_1}{\sqrt{2}}\left(2\Delta v+\Delta q\right)\\
\frac{\delta_2}{\sqrt{2}}\Delta q\\
\end{array}\right]
\]
Since $\Delta v>0$, it is immediate to note that
$
\b{z}\notin\{\b{z}^{*(j)}\}\oplus\mathbb{R}^2_{\succ0}\iff \Delta q<0$,
which means 
\begin{equation}
\label{eq:bounds-proof-1}
\b{z}\notin\{\b{z}^{*(j)}\}\oplus\mathbb{R}^2_{\succ0}\iff q<z_q^{*(j)}+|v-z_v^{*(j)}|
\end{equation}
If $\Delta v<0$, i.e. $z_v^{*(j)}>v$, we have $q=z_q^{*(j)}-\Delta v+\Delta q$. Along the same lines, in this case we obtain:
\[
\b{z}=
\b{z}^{*(j)}+\left[
\begin{array}{l}
\frac{\delta_1}{\sqrt{2}}\Delta q\\
\frac{\delta_2}{\sqrt{2}}\left(2|\Delta v|+\Delta q\right)\\
\end{array}\right]
\]
and  we reach again the statement \eqref{eq:bounds-proof-1}. Now, for $\b{z}$ to be non-dominated by any of the samples $\b{z}^{*(j)}$, the condition \eqref{eq:bounds-proof-1} must hold for all of them, i.e.:
\begin{equation}\label{eq:bounds-proof-2}
\b{z}\notin\bigcup\limits_{\b{z}^{*(j)}\in\mathcal{D}}\{\b{z}^{*(j)}\}\oplus\mathbb{R}^2_{\succ0}\iff q<\min\limits_{\b{z}^{*(j)}\in\mathcal{D}}z_q^{*(j)}+|v-z_v^{*(j)}|
\end{equation}
Note that if $q=\min\limits_{\b{z}^{*(j)}\in\mathcal{D}}z_q^{*(j)}+|v-z_v^{*(j)}|$, then we obtain either $\b{z}=\b{z}^{*(j)}+\left[0,\frac{2\delta_2}{\sqrt{2}}|\Delta v|\right]^T$ or $\b{z}=\b{z}^{*(j)}+\left[\frac{2\delta_1}{\sqrt{2}}|\Delta v|,0\right]^T$ for some $j$: in both cases, the considered point can not belong to the PF. Thus, the function
\begin{equation}\label{eq:bounds-proof-3}
\min\limits_{\b{z}^{*(j)}\in\mathcal{D}}z_q^{*(j)}+|v-z_v^{*(j)}|=\sup_{h\in{{FPFS}_{\mathcal{D}}}}{h(v)}=\bar{h}_{\mathcal{D}}(v)
\end{equation}
is the optimal upper bound of $FPFS_\mathcal{D}$.\\
To find the optimal lower bound, we repeat the reasoning by taking this time $q=z_q^{*(j)}-|v-z_v^{*(j)}|+\Delta q$ for each $j$. Considering 
$v=z_v^{*(j)}+\Delta v,\,\Delta v>0$ and applying the inverse transformation to $[v,q]^T$ we obtain:
\[
\b{z}=
\b{z}^{*(j)}+\left[
\begin{array}{l}
\frac{\delta_1}{\sqrt{2}}\Delta q\\
\frac{\delta_2}{\sqrt{2}}\left(-2\Delta v+\Delta q\right)\\
\end{array}\right]
\]
Since $\Delta v>0$, we note that
$
\b{z}\notin\{\b{z}^{*(j)}\}\oplus\mathbb{R}^2_{\prec0}\iff \Delta q>0$,
which means 
\begin{equation}\nonumber
\b{z}\notin\{\b{z}^{*(j)}\}\oplus\mathbb{R}^2_{\succ0}\iff q>z_q^{*(j)}-|v-z_v^{*(j)}|
\end{equation}
The same conclusions is reached if $\Delta v<0$, in this case we have in fact:
\[
\b{z}=
\b{z}^{*(j)}+\left[
\begin{array}{l}
\frac{\delta_1}{\sqrt{2}}\left(-2|\Delta v|+\Delta q\right)\\
\frac{\delta_2}{\sqrt{2}}\Delta q\\
\end{array}\right]
\]
Thus, considering all samples, we obtain:
\begin{equation}\label{eq:bounds-proof-5}
\b{z}\notin\bigcup\limits_{\b{z}^{*(j)}\in\mathcal{D}}\{\b{z}^{*(j)}\}\oplus\mathbb{R}^2_{\prec0}\iff q>\max\limits_{\b{z}^{*(j)}\in\mathcal{D}}z_q^{*(j)}-|v-z_v^{*(j)}|
\end{equation}
i.e., in order for $\b{z}$ to be non-dominating any of the PF samples in $\mathcal{D}$, the value of $q$ must lie above the following function (optimal lower bound):
\begin{equation}\nonumber
\max\limits_{\b{z}^{*(j)}\in\mathcal{D}}z_q^{*(j)}-|v-z_v^{*(j)}|=\inf_{h\in{{FPFS}_{\mathcal{D}}}}{h(v)}=\underline{h}_{\mathcal{D}}(v)
\end{equation}
What remains to be shown is that the optimal bounds can be evaluated by considering just the samples $\b{z}^{*(j^-)}(v),\,\b{z}^{*(j^+)}(v)$ obtained from \eqref{eq:indexes generating bounds}. Considering together conditions \eqref{eq:bounds-proof-2} and \eqref{eq:bounds-proof-5}, we obtain (see \eqref{eq:inactive-active zones global}):
\begin{equation}\label{eq:bounds-proof-7}
\begin{array}{rcl}
\b{z}\in\mathcal{A}_\mathcal{D}&\iff&
\max\limits_{\b{z}^{*(j)}\in\mathcal{D}}z_q^{*(j)}-|v-z_v^{*(j)}|<q\\
&&<\min\limits_{\b{z}^{*(j)}\in\mathcal{D}}z_q^{*(j)}+|v-z_v^{*(j)}|
\end{array}
\end{equation}
Consider now the sample $\b{z}^{*(j^-)}(v)$. Since it belongs to the PF, it also belongs to $\mathcal{A}_\mathcal{D}$. Thus, according to \eqref{eq:bounds-proof-7}, we have 
\[z^{*(j^-)}_q(v)< z^{*(j)}_q+|z_v^{*(j^-)}-z_v^{*(j)}|,\,\forall \b{z}^{*(j)}\in\mathcal{D}\setminus\{\b{z}^{*(j^-)}(v)\},\]
which implies, for any $j$ such that $z_v^{*(j)}<z_v^{*(j^-)}$,
\[
\begin{array}{l}
z^{*(j^-)}_q(v)+|v-z_v^{*(j^-)}|\\< z^{*(j)}_q+|z_v^{*(j^-)}-z_v^{*(j)}|+|v-z_v^{*(j^-)}|\\=z^{*(j)}_q+|v-z_v^{*(j)}|.
\end{array}
\]
The result above derives from the fact that $z_v^{*(j)}<z_v^{*(j^-)}<v$ by the definition of $\b{z}^{*(j^-)}(v)$ \eqref{eq:indexes generating bounds}. Thus, the upper bound $z^{*(j^-)}_q(v)+|v-z_v^{*(j^-)}|$ is lower than the one contributed by any other sample with $v$-coordinate smaller than $z_v^{*(j^-)}$, which can then be disregarded in the computation of \eqref{eq:bounds-proof-3}. A similar reasoning applies to $\b{z}^{*(j^+)}(v)$, which rules out all samples with larger $v$-coordinate, and to the lower bounds. Thus, the optimal bounds can be obtained by considering $\b{z}^{*(j^-)}(v),\,\b{z}^{*(j^+)}(v)$ only,  proving the Theorem.\hfill$\blacksquare$\\

\textit{Proof of Proposition \ref{prop: max global error}}. 
The non-overlapping intervals $R^{(\ell)},\,\ell=1,\ldots,M-1$ cover the whole domain $[z_v^{*(a1)},\,z_v^{*(a2)}]$ except for the points corresponding to the samples, where however the worst-case error is zero because the upper and lower bounds coincide. Denoting 
$\lambda_{\mathcal{D},max}^{(\ell)}=\max\limits_{v\in R^{(\ell)}}\bar{h}_{\mathcal{D}}(v)-\underline{h}_{\mathcal{D}}(v),$
we can then express the global worst-case error bound as:
\begin{equation}\label{eq:max global error proof 1}
\overline{\lambda}_{\mathcal{D}}=\max\limits_{\ell=1,\ldots,M-1}\lambda_{\mathcal{D},max}^{(\ell)}.
\end{equation}
For each $R^{(\ell)}$, considering Theorem \ref{T:bounds}, we see that the upper and lower bounds for any $v\in R^{(\ell)}$ depend only on the points defining its extremes, i.e. $\b{z}^{*(\ell-)},\,\b{z}^{*(\ell+)}$. The value of $\lambda_{\mathcal{D},max}^{(\ell)}$ can be computed analytically on the basis of the $(v,q)$-coordinates of these points, as follows. First, we compute the coordinates $v',\,v''\in R^{(\ell)}$ such that:
\[
\begin{array}{l}
z_q^{*(\ell-)}+(v'-z_v^{*(\ell-)})=z_q^{*(\ell+)}+(z_v^{*(\ell+)}-v')\\
z_q^{*(\ell-)}-(v''-z_v^{*(\ell-)})=z_q^{*(\ell+)}-(z_v^{*(\ell+)}-v'')\\
\end{array}
\]
i.e., the $v$-values at which the upper (respectively lower) bounds, contributed by the two samples delimiting the interval $R^{(\ell)}$, coincide.
With a few manipulations, one obtains:
\begin{equation}\label{eq:max global error proof 2}
\begin{array}{l}
v'=\dfrac{z_v^{*(\ell+)}+z_v^{*(\ell-)}}{2}+\dfrac{Q^{(\ell)}}{2}\\
v''=\dfrac{z_v^{*(\ell+)}+z_v^{*(\ell-)}}{2}-\dfrac{Q^{(\ell)}}{2}
\end{array}
\end{equation}
Thus, if the $q$-coordinates of $\b{z}^{*(\ell-)},\,\b{z}^{*(\ell+)}$ are the same ($Q^{(\ell)}=0$), $v'$ and $v''$ coincide and are equal to the middle point of $R^{(\ell)}$, otherwise they differ and we can have either $v'>v''$ or $v'<v''$, depending on the sign of $Q^{(\ell)}=z_q^{*(\ell+)}-z_q^{*(\ell-)}$. Assume that $Q^{(\ell)}<0$, which implies $v'<v''$. Then, the interval $R^{(\ell)}$ can be split in three sub-intervals: $I=(z_v^{*(\ell-)},\,v'],\,II=(v',v''),\,III=[v'',\,z_v^{*(\ell+)})]$. In intervals $I$ and $III$, applying \eqref{eq:optimal bounds} one can find that both the upper and lower bounds are determined either by the sample $\b{z}^{*(\ell-)}$ (in $I$) or $\b{z}^{*(\ell+)}$ (in $III$), their maximum distance is attained when $v=v'$ (respectively $v=v''$) and it is equal to (using \eqref{eq:max global error proof 2} and considering that we are assuming $Q^{(\ell)}<0$):
\begin{equation}\label{eq:max global error proof 3}
\lambda_{\mathcal{D},max}^{(\ell)}=V^{(\ell)}-|Q^{(\ell)}|
\end{equation}
In interval $II$, the upper bound is determined by $\b{z}^{*(\ell+)}$, and the lower one by sample $\b{z}^{*(\ell-)}$. Using again \eqref{eq:max global error proof 2} and the assumption $Q^{(\ell)}<0$, one can see that the distance between the bounds is constant in the whole interval and equal to \eqref{eq:max global error proof 3}. The same conclusions are reached if $Q^{(\ell)}>0$, thus the quantity $\lambda_{\mathcal{D},max}^{(\ell)}$ \eqref{eq:max global error proof 3} is indeed the worst-case error bound pertaining to the interval $R^{(\ell)}$. By combining \eqref{eq:max global error proof 1} and \eqref{eq:max global error proof 3}, the result is proven.\hfill$\blacksquare$\\

\textit{Proof of Theorem \ref{th:2}}.
    \noindent $(\Leftarrow)$ For any $\hat{\mathcal{P}}\in FPFS_\mathcal{D}$ and any $\hat{\b{z}}\in\hat{\mathcal{P}}$, we consider the corresponding transformed point
    \[
    \left[\begin{array}{c}
        \hat{z}_v \\
        \hat{z}_q
    \end{array}\right]=W\hat{\b{z}}
    \]
    where $\hat{z}_v\in[z^{*(a1)}_{v},z^{*(a2)}_{v}]$. Consider the realization    
    $\b{z}^{*}=\mathcal{R}_{SM}(\hat{\b{z}})$ \eqref{eq:SM realization}. The corresponding approximation error is  
    \[
    \b{\varepsilon}=\hat{\b{z}}-\b{z}^{*}=\underbrace{\frac{\delta_{1}\delta_{2}}{\sqrt{2}}
            \begin{bmatrix}
            \frac{1}{\delta_2} &  \frac{1}{\delta_2} \\
             -\frac{1}{\delta_1} & \frac{1}{\delta_1}\\ 
            \end{bmatrix}}\limits_{W^{-1}}\left[
    \begin{array}{c}
    0\\
    \hat{z}_q- h(\hat{z}_v)
    \end{array}
    \right]
    \]
    thus,
\begin{equation}\label{eq:th2-necessary-1}
|\varepsilon_1|=\frac{\delta_1}{\sqrt{2}} |\hat{z}_q- h(\hat{z}_v)|;\,
|\varepsilon_2|=\frac{\delta_2}{\sqrt{2}} |\hat{z}_q- h(\hat{z}_v)|  
\end{equation}
According to Theorem \ref{T:bounds}, both $h(\hat{z}_v)$ and $\hat{z}_q$ belong to the interval $(\underline{h}_{\mathcal{D}}(\hat{z}_v),\,\bar{h}_{\mathcal{D}}(\hat{z}_v))$, therefore $|\hat{z}_q- h(\hat{z}_v)|<\bar{h}_{\mathcal{D}}(\hat{z}_v)-\underline{h}_{\mathcal{D}}(\hat{z}_v)$. By \eqref{eq:local-global errors} and the Theorem's necessary condition, we further have $\bar{h}_{\mathcal{D}}(\hat{z}_v)-\underline{h}_{\mathcal{D}}(\hat{z}_v)\leq\bar{\lambda}_\mathcal{D}\leq\sqrt{2}$. Using these inequalities in \eqref{eq:th2-necessary-1} we obtain $|\varepsilon_1|<\delta_1$ and $
|\varepsilon_2|<\delta_2$, thus proving the condition.

$(\Rightarrow)$ We prove the sufficient condition by contradiction. First, note that \eqref{eq:robust requirement} can be reformulated, considering the transformed error coordinates $\b{\varepsilon}_{vq}=W\b{\varepsilon}=[\varepsilon_v,\,\varepsilon_q]^T$, as 
\begin{equation}\label{eq:norm1_proof}
{\lVert\b{\varepsilon}_{vq}\rVert}_1<\sqrt{2}
\end{equation}
where ${\lVert\cdot\rVert}_1$ is the $\ell_1$ vector norm.
This condition can be derived from the inverse transformation $\boldsymbol\varepsilon=W^{-1}\boldsymbol\varepsilon_{vq}$:
\begin{align*}
    &\boldsymbol\varepsilon=\frac{\delta_{1}\delta_{2}}{\sqrt{2}}
            \begin{bmatrix}
            \frac{1}{\delta_2} &  \frac{1}{\delta_2} \\
             -\frac{1}{\delta_1} & \frac{1}{\delta_1}\\ 
            \end{bmatrix}
            \begin{bmatrix}
            \varepsilon_v\\ 
            \varepsilon_q
            \end{bmatrix}
            = 
             \begin{bmatrix}
            \frac{\delta_1(\varepsilon_v+\varepsilon_q)}{\sqrt{2}}\\ 
            \frac{\delta_2(\varepsilon_q-\varepsilon_v)} {\sqrt{2}}
            \end{bmatrix}\\
    & |\varepsilon_1|<\delta_1\Leftrightarrow\frac{\delta_1|\varepsilon_v+\varepsilon_q|}{\sqrt{2}}<\delta_1\Leftrightarrow|\varepsilon_v+\varepsilon_q|<\sqrt{2}\\
    & |\varepsilon_2|<\delta_2\Leftrightarrow\frac{\delta_2|\varepsilon_q-\varepsilon_v|}{\sqrt{2}}<\delta_2\Leftrightarrow|\varepsilon_q-\varepsilon_v|<\sqrt{2}\\
    & |\varepsilon_q\pm\varepsilon_v|<\sqrt{2}\Leftrightarrow\lVert\boldsymbol\varepsilon_{vq}\rVert_1<\sqrt{2}
\end{align*}
Suppose, for the sake of contradiction, that $\bar{\lambda}_\mathcal{D}>\sqrt{2}$. Then, we can select $\hat{\mathcal{P}},\mathcal{P}\in FPFS_\mathcal{D}$ and two points $\hat{\b{z}}\in\hat{\mathcal{P}},\,\b{z}^*\in\mathcal{P}$, such that
\begin{equation}\label{eq:th2-sufficient-1}
       \hat{z}_v=z^*_v
       \vee|\hat{z}_q-z^*_q|>\sqrt{2}
\end{equation}
Assume now that there is still a realization $\mathcal{R}$ such that property \eqref{eq:robust requirement} holds, and denote $\tilde{{\b{z}}}^{*}=\mathcal{R}(\hat{\b{z}})$. Further denote $\b{r}=\tilde{{\b{z}}}^{*}-\b{z}^*$ and $[r_v,\,r_q]^T=W\b{r}$. According to \eqref{eq:norm1_proof}, it shall hold that
$|\hat{z}_v-\tilde{z}^*_v| + |\hat{z}_q-\tilde{z}^*_q| < \sqrt{2}$.
Then, using \eqref{eq:th2-sufficient-1} we obtain:
\begin{align*}
& {\lVert\varepsilon_{vq}\rVert}_1=|\hat{z}_v-\tilde{z}^*_v| + |\hat{z}_q-\tilde{z}^*_q| < \sqrt{2} \\
&|\hat{z}_v-z^*_v-r_v| + |\hat{z}_q-z^*_q-r_q| < \sqrt{2}\\
&|r_v| + \sqrt{2}-|r_q| < \sqrt{2}\\
& |r_v|-|r_q| < 0 \Rightarrow |r_v|<|r_q|\Rightarrow\\ &\Rightarrow |\tilde{z}^*_v-z^*_v|<|\tilde{z}^*_q-z^*_q|
\end{align*}
The inequality above
implies:
\[
\tilde{z}^*_q>z^*_q+|\tilde{z}^*_v-z^*_v|\, \lor\,
\tilde{z}^*_q< z^*_q-|\tilde{z}^*_v-z^*_v|
\]
which, in either case, shows that the point $\tilde{z}^*$ is not admissible, since it violates the Lipschitz continuity property of the transformed PF (see Remark \ref{rem:Lipschitz}). In the non-rotated coordinates, it means that $\tilde{z}^*$ either dominates or is dominated by $\b{z}^*$, thus completing the contradictory argument and showing that the condition $\bar{\lambda}_\mathcal{D}\leq\sqrt{2}$ is necessary for the existence of a realization that satisfies \eqref{eq:robust requirement} $\forall\hat{\mathcal{P}},\mathcal{P}\in FPFS_\mathcal{D}$ and, conversely, the latter condition is sufficient to have $\bar{\lambda}_\mathcal{D}\leq\sqrt{2}$.\hfill$\blacksquare$\\ 

\textit{Proof of Proposition \ref{prop:minimum M}}. 
    From Proposition \ref{prop: max global error}, we need to have $V^{(\ell)}-|Q^{(\ell)}|\leq\sqrt{2},\,\ell=1,\ldots,\,M-1$. From the same result, we can see that the worst case with respect to the actual values of the samples is when $Q^{(\ell)}=0,\,\forall \ell$, i.e. the rotated PF samples have all the same $q$-coordinate. Then, considering that the length of the whole domain is $V_a=z^{*(a2)}_v-z^{*(a1)}_v$, the condition \eqref{eq:main result} is satisfied if $M=\ceil*{V_a/\sqrt{2}}+1$ using evenly distributed points in the segment $[z^{*(a2)}_v,z^{*(a1)}_v]$ including the extremes. Note that any other allocation of samples would result in a larger worst-case error. Finally, we re-write $V_a$ as a function of $\Delta_1,\,\Delta_2,\,\delta_1,\,\delta_2$ by applying the transformation $W$ to the anchor points:
    \begin{equation}\nonumber
    \begin{split}
            &V_a=z^{*(a2)}_v-z^{*(a1)}_v\\
            &V_a=\frac{\sqrt{2}}{2}\left(\frac{z_1^{*(a2)}}{\delta_1}-\frac{z_2^{*(a2)}}{\delta_2}\right)-\frac{\sqrt{2}}{2}\left(\frac{z_1^{*(a1)}}{\delta_1}-\frac{z_2^{*(a1)}}{\delta_2}\right)\\
            &\frac{V_a}{\sqrt{2}}=\frac{1}{2}\left(\frac{\Delta_1}{\delta_1}+\frac{\Delta_2}{\delta_2}\right)
    \end{split}
\end{equation}
\hfill$\blacksquare$\\

\textit{Proof of Proposition \ref{prop:central estimate}}. 
The first part of the result can be demonstrated as a corollary of Theorem \ref{th:2}: since the worst-case error obtained by the central approximation is  $\frac{1}{2}\overline{\lambda}_{\mathcal{D}}$, applying the same reasoning of the proof of Theorem \ref{th:2} one can conclude that $\overline{\lambda}_{\mathcal{D}}\leq2\sqrt{2}$ is necessary and sufficient for the robust error guarantees. 
    The second part, pertaining to the number of samples, can be demonstrated with the same arguments as Proposition \ref{prop:minimum M}, considering that the value of $V^{(\ell)}$ shall now be smaller than $2\sqrt{2}$ for all $\ell=1,\ldots,M-1$.\hfill$\blacksquare$    

\bibliographystyle{ieeetr}
\bibliography{ROBBO_literature}

\begin{thebibliography}{10}

\bibitem{BOP1}
L.~Tang, Y.~Li, D.~Bai, T.~Liu, and L.~C. Coelho, ``Bi-objective optimization for a multi-period {COVID}-19 vaccination planning problem,'' {\em Omega}, vol.~110, p.~102617, 2022.

\bibitem{BOP3}
E.~Demir, T.~Bektaş, and G.~Laporte, ``The bi-objective pollution-routing problem,'' {\em European Journal of Operational Research}, vol.~232, no.~3, pp.~464--478, 2014.

\bibitem{BOP4}
S.~Wang, M.~Liu, F.~Chu, and C.~Chu, ``Bi-objective optimization of a single machine batch scheduling problem with energy cost consideration,'' {\em Journal of Cleaner Production}, vol.~137, pp.~1205--1215, 2016.

\bibitem{MultiobjectiveMetaheuristicOptimization2020a}
A.~Rodríguez-Molina, E.~Mezura-Montes, M.~G. Villarreal-Cervantes, and M.~Aldape-Pérez, ``Multi-objective meta-heuristic optimization in intelligent control: {A} survey on the controller tuning problem,'' {\em Applied Soft Computing}, vol.~93, p.~106342, Aug. 2020.

\bibitem{marlerSurveyMultiobjectiveOptimization2004}
R.~Marler and J.~Arora, ``Survey of multi-objective optimization methods for engineering,'' {\em Structural and Multidisciplinary Optimization}, vol.~26, no.~6, pp.~369--395, 2004.

\bibitem{dasCloserLookDrawbacks1997}
I.~Das and J.~E. Dennis, ``A closer look at drawbacks of minimizing weighted sums of objectives for pareto set generation in multicriteria optimization problems,'' {\em Structural optimization}, vol.~14, no.~1, pp.~63--69, 1997.

\bibitem{dasNormalBoundaryIntersectionNew1998}
I.~Das and J.~E. Dennis, ``Normal-boundary intersection: A new method for generating the pareto surface in nonlinear multicriteria optimization problems,'' {\em {SIAM} Journal on Optimization}, vol.~8, no.~3, pp.~631--657, 1998.
\newblock Publisher: Society for Industrial and Applied Mathematics.

\bibitem{ismail-yahayaEffectiveGenerationPareto2002}
A.~Ismail-Yahaya and A.~Messac, ``Effective generation of the {P}areto frontier using the normal constraint method,'' in {\em 40th {AIAA} Aerospace Sciences Meeting \& Exhibit}, American Institute of Aeronautics and Astronautics, 2002.

\bibitem{CHIANDUSSI2012912}
G.~Chiandussi, M.~Codegone, S.~Ferrero, and F.~Varesio, ``Comparison of multi-objective optimization methodologies for engineering applications,'' {\em Computers \& Mathematics with Applications}, vol.~63, no.~5, pp.~912--942, 2012.

\bibitem{alianofilhoExactScalarizationMethod2021}
A.~Aliano~Filho, A.~C. Moretti, M.~V. Pato, and W.~A. de~Oliveira, ``An exact scalarization method with multiple reference points for bi-objective integer linear optimization problems,'' {\em Annals of Operations Research}, vol.~296, no.~1, pp.~35--69, 2021.

\bibitem{emmerichTutorialMultiobjectiveOptimization2018}
M.~T.~M. Emmerich and A.~H. Deutz, ``A tutorial on multiobjective optimization: fundamentals and evolutionary methods,'' {\em Natural Computing}, vol.~17, no.~3, pp.~585--609, 2018.

\bibitem{pereiraReviewMultiobjectiveOptimization2022a}
J.~L.~J. Pereira, G.~A. Oliver, M.~B. Francisco, S.~S. Cunha, and G.~F. Gomes, ``A review of multi-objective optimization: Methods and algorithms in mechanical engineering problems,'' {\em Archives of Computational Methods in Engineering}, vol.~29, no.~4, pp.~2285--2308, 2022.

\bibitem{cohonGeneratingMultiobjectiveTradeoffs1979}
J.~L. Cohon, R.~L. Church, and D.~P. Sheer, ``Generating multiobjective trade-offs: An algorithm for bicriterion problems,'' {\em Water Resources Research}, vol.~15, no.~5, pp.~1001--1010, 1979.

\bibitem{kimAdaptiveWeightedSum2006}
I.~Y. Kim and O.~L. de~Weck, ``Adaptive weighted sum method for multiobjective optimization: a new method for {P}areto front generation,'' {\em Structural and Multidisciplinary Optimization}, vol.~31, no.~2, pp.~105--116, 2006.

\bibitem{klamrothUnbiasedApproximationMulticriteria2003}
K.~Klamroth, J.~Tind, and M.~M. Wiecek, ``Unbiased approximation in multicriteria optimization,'' {\em Mathematical Methods of Operations Research ({ZOR})}, vol.~56, no.~3, pp.~413--437, 2003.

\bibitem{bokrantzAlgorithmApproximatingConvex2013a}
R.~Bokrantz and A.~Forsgren, ``An algorithm for approximating convex pareto surfaces based on dual techniques,'' {\em {INFORMS} Journal on Computing}, vol.~25, no.~2, pp.~377--393, 2013.

\bibitem{ruzikaApproximationMethodsMultiobjective2005}
S.~Ruzika and M.~M. Wiecek, ``Approximation methods in multiobjective programming,'' {\em Journal of Optimization Theory and Applications}, vol.~126, no.~3, pp.~473--501, 2005.

\bibitem{eichfelderAdvancementsComputationEnclosures2023b}
G.~Eichfelder and L.~Warnow, ``Advancements in the computation of enclosures for multi-objective optimization problems,'' {\em European Journal of Operational Research}, vol.~310, no.~1, pp.~315--327, 2023.

\bibitem{eichfelderApproximationAlgorithmMultiobjective2022}
G.~Eichfelder and L.~Warnow, ``An approximation algorithm for multi-objective optimization problems using a box-coverage,'' {\em Journal of Global Optimization}, vol.~83, no.~2, pp.~329--357, 2022.

\bibitem{payneInteractiveRectangleElimination1980}
A.~Payne and E.~Polak, ``An interactive rectangle elimination method for biobjective decision making,'' {\em {IEEE} Transactions on Automatic Control}, vol.~25, no.~3, pp.~421--432, 1980.

\bibitem{payne1993efficient}
A.~N. Payne, ``Efficient approximate representation of bi-objective tradeoff sets,'' {\em Journal of the Franklin Institute}, vol.~330, no.~6, pp.~1219--1233, 1993.

\bibitem{campigottoActiveLearningPareto2014}
P.~Campigotto, A.~Passerini, and R.~Battiti, ``Active learning of pareto fronts,'' {\em {IEEE} Transactions on Neural Networks and Learning Systems}, vol.~25, no.~3, pp.~506--519, 2014.

\bibitem{binoisQuantifyingUncertaintyPareto2015}
M.~Binois, D.~Ginsbourger, and O.~Roustant, ``Quantifying uncertainty on pareto fronts with gaussian process conditional simulations,'' {\em European Journal of Operational Research}, vol.~243, no.~2, pp.~386--394, 2015.

\bibitem{tianLocalModelBasedPareto2023a}
Y.~Tian, L.~Si, X.~Zhang, K.~C. Tan, and Y.~Jin, ``Local {Model}-{Based} {Pareto} {Front} {Estimation} for {Multiobjective} {Optimization},'' {\em IEEE Transactions on Systems, Man, and Cybernetics: Systems}, vol.~53, pp.~623--634, Jan. 2023.

\bibitem{emmerichConeBasedHypervolumeIndicators2013}
M.~Emmerich, A.~Deutz, J.~Kruisselbrink, and P.~K. Shukla, ``Cone-based hypervolume indicators: Construction, properties, and efficient computation,'' in {\em Evolutionary Multi-Criterion Optimization}, Lecture Notes in Computer Science, pp.~111--127, Springer, 2013.

\bibitem{nakayamaBasicConceptsMultiobjective2009}
H.~Nakayama, Y.~Yun, and M.~Yoon, ``Basic concepts of multi-objective optimization,'' in {\em Sequential Approximate Multiobjective Optimization Using Computational Intelligence}, Vector Optimization, pp.~1--15, Springer, 2009.

\bibitem{logistFastParetoSet2010a}
F.~Logist, B.~Houska, M.~Diehl, and J.~Van~Impe, ``Fast {Pareto} set generation for nonlinear optimal control problems with multiple objectives,'' {\em Structural and Multidisciplinary Optimization}, vol.~42, pp.~591--603, Oct. 2010.

\bibitem{MILANESE2004957}
M.~Milanese and C.~Novara, ``Set membership identification of nonlinear systems,'' {\em Automatica}, vol.~40, no.~6, pp.~957--975, 2004.

\bibitem{TRAUB80}
J.~Traub and J.~F. Woźniakowski, {\em A General Theory of Optimal Algorithms}.
\newblock Academic Press, 1980.

\bibitem{burachikNewScalarizationTechnique2014a}
R.~S. Burachik, C.~Y. Kaya, and M.~M. Rizvi, ``A new scalarization technique to approximate pareto fronts of problems with disconnected feasible sets,'' {\em Journal of Optimization Theory and Applications}, vol.~162, no.~2, pp.~428--446, 2014.

\bibitem{messacNormalizedNormalConstraint2003}
A.~Messac, A.~Ismail-Yahaya, and C.~Mattson, ``The normalized normal constraint method for generating the pareto frontier,'' {\em Structural and Multidisciplinary Optimization}, vol.~25, no.~2, pp.~86--98, 2003.

\bibitem{ROBBOweb}
R.~Boffadossi, ``Robbo demo web page,'' 2025.
\newblock Available at \url{https://www.sas-lab.deib.polimi.it/?p=1296}.

\bibitem{MAGNI20011351}
L.~Magni, G.~Nicolao, L.~Magnani, and R.~Scattolini, ``A stabilizing model-based predictive control algorithm for nonlinear systems,'' {\em Automatica}, vol.~37, no.~9, pp.~1351--1362, 2001.

\end{thebibliography}

\end{document}